\RequirePackage{fix-cm}
\documentclass{svjour3}                     
\smartqed  
\usepackage{amsmath}
\usepackage{amssymb}
\usepackage{latexsym}
\usepackage{color}
\usepackage{graphicx}
\newcommand{\pr}{\mathbb{P}}

\newcommand{\E}{\mathbb{E}}
\newcommand{\ep}{\hbox{ }\hfill{ ${\cal t}$~\hspace{-5.1mm}~${\cal u}$   } }
\newcommand{\eps}{\varepsilon}

\begin{document}

{\title{Reducing the debt: Is it optimal to outsource an investment?}}

\titlerunning{Reducing the debt}        

\author{G.E.~Espinosa  \and ~C.~Hillairet  \and ~B.~Jourdain   \and~M.~Pontier}
%

\institute{Universit\'e Paris-Est, CERMICS (ENPC), \email{gilles-edouard.espinosa@polytechnique.org} \at
  \\ 
 \and
        CMAP, Ecole Polytechnique,  \email{caroline.hillairet@polytechnique.edu} \at This research benefited from the support of the Chaires ``D\'eriv\'es
du Futur'' and  ``Finance et D\'eveloppement durable''.
         \\
\and
          Universit\'e Paris-Est, CERMICS (ENPC), INRIA,  \email{jourdain@cermics.enpc.fr} \at This research benefited from the support of the Eurostars E!5144-TFM project and of the Chaire ``Risques Financiers``, Fondation du risque.
 \\
          \and
          IMT  Universit\'e Paul Sabatier,   \email{pontier@math.univ-toulouse.fr} 
         \\   
}
    

%

\date{Received: date / Accepted: date}
\maketitle

\abstract
{This paper deals with the problem of outsourcing the debt for a big investment, according two 
situations: either the firm outsources both the investment (and the associated debt) and the exploitation to another firm (for example a  private consortium),
or the firm supports the debt and the investment but outsources the exploitation.
We prove the existence of Stackelberg and Nash  equilibria between the firms,
 in both situations. We compare
the benefits of these  contracts, theorically and numerically. We conclude with a study of what 
happens in case of incomplete information, in the sense that the risk aversion coefficient of each partner may be unknown by the other partner.
}\\


\section{Introduction}

With the significant increase in recent years of public debt in many developed countries, together with the associated concerns related to possible defaults of some of them, the question of financing public projects is more than ever at the center of economic and political considerations. To overcome this problem, leveraging on the private sector appears at first glance as a good idea. This type of Public-Private Partnership (PPP) was initiated in the United Kingdom in 1992, under the name Private Financing Initiative (PFI), and has widely been used since then, so that it represented one third of all public investments made in the UK during the period 2001-2006. It has also been used in many other countries, in particular in Europe, Canada and in the United States, to finance hospitals, prisons or stadiums among others. It is also recommended by the OECD. We refer among many other references to \cite{EFG} for an overview of the extent of PPPs in Europe and in the US, to the website of the European PPP Expertise Center (EPEC) or the website of the National Council for Public-Private Partnerships, and to \cite{OECD} for a global overview made by the OECD.

However, as emphasized by the recent discussions in the UK, although the benefits of this type of partnership are mainly admitted, there are still many concerns about its drawbacks (see \cite{WBI} for a detailed overview). Even though some drawbacks are of political, social or behavioral natures, others are purely economic and are the ones that we are interested in. More precisely, we would like to answer the following question: from an economic point of view, and taking into account the constraints that a country faces when issuing a new amount of debt, is it optimal for this country to finance a public project via a private investment?

Although already of a big interest, this question does not need to be restricted to debt issuance by a country but can be generalized to any economic agent, be it a country or a firm. Indeed, any firm has some constraints on its debt level for several reasons. In some cases, banks will simply not allow a company to borrow enough money to sustain a very expensive project. But even if it is not the case, since the debt level appears on the company's accounts, issuing too much debt will affect the opinion and confidence of investors, and in particular deteriorate its rating. This can lead to a higher credit spread when issuing new bonds, difficulties to increase the capital of the firm, a drop of the company's stock price, dissatisfaction of shareholders or in the worst case, bankruptcy. We can cite some concrete examples where the dilemma between investing directly or resorting to outside investment can occur: owning or renting offices or factories, owning or leasing trucks, trains or planes, some industrial machines or some office materials (such as computers).

Therefore we will consider in this paper the problem of outsourcing from the debt point of view. Since the question of outsourcing some operations has already been widely studied and our aim is only to study the relevance of outsourcing an investment in order to reduce the debt of a firm (or economic agent), we will compare two situations where the operations are always outsourced. In the first one, the firm outsources both the investment/debt and the operations, while in the second one, the firm supports the debt and the investment but outsources the operations. In both situations, the investment is covered by issuing a debt at time 0 but the cost of borrowing of the outsourcer and the outsourcee  may be different. More specifically, we suppose that the outsourcer faces some important constraints if he has to issue a new debt, stylized by  a convex function $f$ modeling his aversion for debt issuance. That is why he considers the possibility  to outsourcing the investment to a firm which has less constraints.  

In  \cite{IMP}, Iossa, Martimort and Pouyet give some results on the comparison of the costs and benefits associated to PPP. Hillairet and Pontier \cite{CaroMo} propose a study on PPP  and their relevance,  assuming the eventuality of a default of the counterparty, but  they do not take into account the government debt aversion. However,  the attractivity for government of PPP contracts relies obviously on  the short term opportunity gain to record infrastructure assets out of the government's book.  To our knowledge, there   does not exist any
reference in Mathematics area.   In Economics, a narrow strand of literature is  dedicated to the discussion of Build-Operate-Transfer (BOT) concession  contracts, which is a frequent form of PPP.  Under BOT contracts the private sector builds and operates an infrastructure project for a well defined concession period and then transfers it to public authorities. The attractiveness of  BOT contracts to governments stems from the possibility to limit governmental spending by shifting the investment costs to a private consortium.  
In \cite{AuPi1}, Auriol and Picard discuss the choice of BOT contracts when governments and consortia do not share the same information about the cost parameter during the project life. They summarize the government's financial constraint by its ``shadow'' cost of public funds, which reflects the macro-economic constraints that are imposed on national governments' surplues and debt levels by supranational institutions such as the I.M.F. Using linear demand functions and uniform cost distributions, they compute theoretical values of shadow costs that would entice governments to choose BOT concessions contracts.  
Our approach is different from the modelization and the resolution point of view. 

This paper studies, in a two-period setting, two kinds of equilibria  between risk averse firms. The first one is a Stackelberg equilibrium with the outsourcer as leader, which corresponds to  the more classic setting for outsourcing situations. The second one is  a Nash equilibrium. It may correspond to an outsourcing between two economic entities within a same group: in this case, the two  entities make their decisions simultaneously and a Nash equilibrium is more favourable for the group than a Stackelberg one. For both situations where the investment is outsourced or not, Stackelberg and Nash equilibria are characterized. A theorical comparison is provided for the Nash equilibrium, from the point of view of the outsourcer: we check  that the investment should be outsourced if the outsourcee has a lower cost of capital or if the outsourcer has a high debt aversion.  For the Stackelberg equilibrium, the comparison is done numerically on a concrete example. The analysis is extended to an incomplete information setting in which the firms do not have perfect knowledge of the preference of their counterparty. To model the social need of the investment, the outsourcer gets a penalty if the outsourcee does not accept the contract.  \\

The present paper is organised as follows. In Section \ref{pb} we set the problem of outsourcing between two firms and we define the optimization problems in Situation 1,  in which the firm outsources both the investment/debt and the operations, and in Situation 2 , in which the firm supports the debt and the investment but outsources the operations.
Section \ref{sec:result} provides the main results concerning Nash and Stackelberg equilibria, the comparison between the two situations, and the case of incomplete information. The proofs are postponed in Appendices.
 We provide in Section \ref{sec:example}  a numerical example to better investigate the quantitative effects of the model. 
Appendix \ref{sec:proof1} is devoted to the proofs  of existence and characterization of Nash and Stackelberg equilibria in  Situation 1, then Appendix \ref{sec:proof2} does the same  in Situation 2. Appendix \ref{sec:proof3} concerns the proofs of the 
comparison results between the two situations and Appendix \ref{proofIncInf} the results obtained in incomplete information.

\section{Problem formulation}
\label{pb}


\subsection{Costs and revenue}

Consider two firms. Firm $I$ is the one who wants to reduce its debt and therefore considers the possibility of outsourcing an investment to a second firm $J$. In any case, firm $J$ is the one that will support the operational cost of the project, on the time horizon $T$. Let  the operational cost on the time-interval $[0,T]$ $C^{op}$ be given by:
 \begin{equation}
 \label{EqOperationalCost}
 C^{op}=\mu-\varphi(e)-\delta \psi(a),
 \end{equation}
where
\begin{itemize}
\item $\mu$ is the ``business as usual'' cost, such that $\mathbb{E}[ \mu]$ represents the ``average'' benchmark cost (it takes into account the price of commodities, employees, rents...). We assume that $\mu$ is a non constant random variable bounded from below by a finite constant $\underline{\mu}$ on a probability space $(\Omega, F_T,\pr)$
such that
\begin{equation}
\label{hypmu}
\forall \lambda\in\mathbb{R},\;\mathbb{E}\left[e^{\lambda \mu}\right]<+\infty.
\end{equation}
Notice that this implies that the function $\lambda\mapsto \mathbb{E}\left[e^{\lambda \mu}\right]$ is infinitely differentiable.

\item $e$ is a non-negative parameter and represents the effort made   on the time-interval $[0,T]$ in order to reduce the operational cost such as logistics improvements, research and development, maintenance or more efficient or less workers. $e$  will in general have a social impact for firm $I$,
\item $\delta\in\mathbb{R}$ represents the impact of the quality of the investment on the reduction of operational costs. We do not impose any restriction on the sign of $\delta,$ since, as suggested in \cite{IMP}, both signs can make sense depending on the situation. Indeed, when constructing a building, using more expensive material usually brings less maintenance costs and therefore a positive delta. On the contrary, for a hospital, using more sophisticated (and expensive) machines can bring bigger maintenance costs and a negative delta.
\item 
$a\geq 0$ is the effort done at time $0$ to improve  the (initial) quality of the investment, improving in the meantime the social value of the project. 
 Depending on $\delta$, $a$ affects positively, negatively or does not affect the operational cost. Remark is that $a$ has the same dimension as the effort $e$,
\item $\varphi:\mathbb{R}_+\to\mathbb{R}_+$ and $\psi:\mathbb{R}_+\to\mathbb{R}_+$, modeling the respective impacts of the efforts $e$ and $a$ respectively on the operational cost $C^{op}$, 
 are $C^1$, increasing and strictly concave functions, satisfying the Inada conditions $\varphi'(0)=+\infty$ and $\varphi'(\infty)=0$, $\psi'(0)=+\infty$ and $\psi'(\infty)=0$. 
We also assume that $\varphi(\infty)+\delta^+\psi(\infty)< \underline{\mu}$ where $\underline{\mu}$ is the lower bound of the random variable $\mu$ and $\delta^+=\max(\delta,0)$, which ensures that $C^{op} > 0$; as a consequence, 
$\forall (x,y)\in \mathbb{R}_+^2,$ $\varphi(x)+\delta^+\psi(y)< \underline{\mu}.$
 We assume furthermore that  $(\varphi')^{-1}(x) \sim(\psi')^{-1}(x) $ for
  $x \rightarrow 0$
 for technical reasons and to make the computations lighter even if we could relax this last assumption.

\end{itemize}
\noindent
The parameter $e$ is a control for firm $J$, while $a$ is a control for the firm that supports the investment, $I$ or $J$ depending on the situation;  $\mu$ represents  the cost  on the time-interval $[0,T]$.
 \\
The minimal investment required by the project is $D>0$ and if initial effort are done (i.e. if $a >0$)  the total investment is the sum $D+a$. This investment is assumed to be entirely covered by issuing a debt with horizon $T$ at time $0$. To take into account the possibility that the cost of borrowing is in general not the same for different firms, we denote the respective non-negative constant interest rates of firms $I$ and $J$ by $r_I$ and $r_J$.
 On the time-interval $[0,T]$, the amount to be reimbursed  by the borrower $K\in\{I,J\}$ is $(1+r_K)(D+a)$.
  \\
Finally, we need to add the remaining costs on $[0,T]$ coming from the effort $e$ as well as the maintenance costs denoted by $m$ :
  \begin{equation}
  \label{EqOtherCost}
 C^m=e+m.
 \end{equation}
The maintenance cost $m$ is a non-negative parameter and will have a social impact for firm $I$. 

\vspace{5mm}

Since firm $I$ gives to firm $J$ either a rent or the right to exploit the project on $[0,T]$, in both cases we can consider a random variable $R$ which corresponds to the endowment for firm $J$ and the rent or shortfall for firm $I$, on $[0,T]$. This variable is computed using a reasonably simple rule, decided at  $t=0$ and subject to a control of firm $I$. In reality, in such contracts, the endowment can be indexed on the price of commodities in the case of transportation or on a real-estate index for the rent of a building. Since firm $I$ wants to have a project of good quality as well as a well maintained project, we assume that $R$ is a non-negative random variable and depends on both $C^{op}$ and the maintenance cost in the following way:
  \begin{equation}
  \label{EqEndowment}
 R=\alpha+\beta C^{op}+\gamma g(m),
 \end{equation}
with $\alpha\geq 0$, $\beta\in\mathbb{R}$, $\gamma\geq 0$ and $g$ is a $C^1$, increasing and strictly concave function on $\mathbb{R}_+$, such that $g'(0)=\infty$ and $g'(\infty)=0$. Moreover, we assume that
 \begin{equation}
 \label{hypginfty}
 m_0:=\inf\{m>0:g(m)>0\}<+\infty.
 \end{equation}
The constants $\alpha$, $\beta$ and $\gamma$ are controls of firm $I$.
We do not put any randomness in the coefficients $\alpha$, $\beta$ and $\gamma$ of $R$ since we consider that they are defined at time $t=0$ by a contract between firms $I$ and $J$. All the randomness in $R$ comes from the operational cost term $C^{op}$. Still, this model allows for an indexation on a benchmark such as the price of commodities or inflation through this dependence on operational costs.

\subsection{Optimization problems}

We assume that the risk aversions of firm $I$ and $J$ are represented respectively by the
 exponential utility functions $U(x)=-e^{-ux}$ and $V(x)=-e^{-vx}$, $x\in\mathbb{R}$, 
 with $u,v>0$.
\\
We consider two different situations: in Situation 1, firm $J$ supports the debt and
 takes care of the exploitation; its controls are $a$, $e$ and $m$, whereas the controls 
 of firm $I$ are $\alpha$, $\beta$ and $\gamma$. In Situation 2, firm $J$ only takes care 
 of the exploitation, its controls are $e$ and $m$, whereas the controls of firm $I$ are
  $a$, $\alpha$, $\beta$ and $\gamma$. Firm $I$ is the one that chooses between the two situations. The optimization problems for firm $J$ respectively in Situation 1 and 2 
  are $\sup_{(a,e,m)}J^1(a,e,m)$ and $\sup_{(e,m)}J^2(e,m)$ respectively, where:
\begin{eqnarray}
\label{defJv1}
 J^1(a,e,m)&=&\mathbb{E} \left[V\big(\alpha+(\beta-1)C^{op}-e+\gamma g(m)-m-(1+r_J)(D+a)\big)\right]
\\
 J^2(e,m)&=&\mathbb{E}\left[V\big(\alpha+(\beta-1)
C^{op}-e+\gamma g(m)-m\big)\right]
\label{defJv2}
\end{eqnarray}
recalling $C^{op}=\mu-\varphi(e)-\delta \psi(a)$.

In these optimization problems, we have assumed that the controls of firm $I$ are given (they have no reason to be the same in the two cases). We look for controls in the following admissible sets: $e,$ $m$ and $a$ are  non-negative constants. The eventuality  that firm $J$ does not accept the contract will be  taken into account  in the constraints of the optimization problem for firm $I$. 

\vspace{5mm}

\noindent On the other hand, we consider that the project has an initial ``social'' value $b^a(a)$   for firm $I$, and  a good maintenance also represents a social benefit $b^m(m)$. The benefits of the efforts on operational costs are modelled through the  function $b^e$. We also introduce a penalization function $f$ representing the aversion for debt issuance of firm $I$ (firm $J$ has no debt aversion). Those functions   satisfy the following hypotheses
\begin{itemize}
\item $b^a:\mathbb{R}_+\to\mathbb{R}_+$ is a $C^1$, increasing and concave function.  $(b^a)'(0)>0$, possibly infinite, $(b^a)'(\infty)=0$ and  $b^a(\infty)<\infty$.
\item $b^m:\mathbb{R}_+\to\mathbb{R}_+$ is  a $C^1$, increasing and concave function, such that  $(b^m)'(0)=\infty$ and $(b^m)'(\infty)=0$.
\item $b^e:\mathbb{R}_+\to\mathbb{R}_+$ is  a $C^1$, increasing and concave function,  such that $(b^e)'(0)=\infty$ and $(b^e)'(\infty)=0$
\item $f$ is an increasing and strictly convex function, satisfying  $f'(\infty)=\infty$.
\end{itemize}

Therefore we write the optimization problem for firm $I$ in both situations as
$\sup_{(\alpha,\beta,\gamma)}I^1(\alpha,\beta,\gamma)$ and $\sup_{(a,\alpha,\beta,\gamma)}I^2(a,\alpha,\beta,\gamma)$ where:
\begin{eqnarray*}
I^1(\alpha,\beta,\gamma)&=&\mathbb{E}\left[b^a(a)+
U\big(b^m(m)+b^e(e)-\alpha-\beta C^{op}-\gamma g(m)\big)\right]
\\
I^2(a,\alpha,\beta,\gamma)&=&\mathbb{E}\left[b^a(a)-f\big((1+r_I)(D+a)\big)
+U\big(b^m(m)+b^e(e)-\alpha-\beta
C^{op}-\gamma g(m)\big)\right].
\end{eqnarray*}
\noindent
Hypotheses on $b^a$ and $f$ imply that $F(a):=b^a(a)-f((1+r_I)(D+a))$
is strictly concave,  satisfies $F'(\infty)=-\infty$ and $F(\infty)=-\infty$.
Finally we assume that  $F'(0)>0$, possibly infinite.
The admissible sets are:
\\
-  in Situation 1, $\alpha\geq 0$, $\beta\in\mathbb{R}$, $\gamma\geq 0$ and such that
{\small\begin{equation}
\hspace{-0.5cm} \label{admiset1}
{\mathbb{E} \left[V\big(\alpha+(\beta-1)(\mu-\varphi(e)-\delta \psi(a))- e+\gamma g(m)-m-(1+r_J)(D+a)\big)\right]\geq V(0)};
\end{equation}
}
- in Situation  2, $a\geq 0$, $\alpha\geq 0$, $\beta\in\mathbb{R}$, $\gamma\geq 0$ and such that
\begin{equation}\label{admiset2}
\mathbb{E}\left[V\big(\alpha+(\beta-1)(\mu-\varphi(e)-\delta \psi(a))-e+\gamma g(m)-m\big)\right]\geq V(0).
\end{equation}
The constraint ensures that firm $J$ will accept the contract, since it is better or equal for it than doing nothing.

\vspace{5mm}

The most natural equilibrium to be considered in outsourcing situations is a Stackelberg equibrium with firm $I$ as leader. It corresponds for example  to a government which outsources the construction of a stadium, or to an industrial group which wants to outsource its trucks/trains to a  transport company. Nevertheless, within a group,  a given entity may be interested in outsourcing its debt to another entity of the same group. In this situation,  a Nash equilibrium is more relevant.  Therefore  this paper addresses those two kinds of equilibria.  




\begin{remark}
A Stackelberg equilibrium with firm $J$ leader is never relevant from an economical perspective, since it is never the outsourcee who makes an offer to the outsourcer. It is also not interesting from a mathematical point of view. Indeed, 
in Situation $i\in\{1,2\}$, since $C^{op}=\mu-\varphi(e)-\delta\psi(a)>0$, $I^{i}$ is decreasing with respect to $\beta$, while $J^{i}$ is increasing with respect to $\beta$. Therefore if firm $J$ is the leader, for any choice of its controls, firm $I$'s optimal controls will always bind the constraint $J^{i}\geq V(0)$. 
\end{remark}

\section{Main results}\label{sec:result}
The best responses of firm $J$ to given controls of firm $I$ turn out to be easily derived. That is why we first present them, before stating our main results concerning Nash and Stackelberg equilibria where these best responses appear. The proofs of the main results are postponed in Appendices \ref{sec:proof1}, \ref{sec:proof2}, \ref{sec:proof3}, \ref{proofIncInf}.

\subsection{Best responses of firm $J$ in Situations $1$ and $2$}\label{bestrespj}
Let us first consider Situation $1$ and suppose that $(\alpha,\beta,\gamma)$ is given in $\mathbb{R}_+\times\mathbb{R}\times \mathbb{R}_+$. For firm $J$ the optimization problem is $\omega$ by $\omega$, and since $U$ is increasing it writes:
$$
\sup_{e\geq 0}\{(1-\beta) \varphi(e)-e\}+\sup_{m\geq 0}\{\gamma g(m)-m\}+\sup_{a \geq 0}\{\delta (1-\beta) \psi(a) -(1+r_J)a \} .
$$
Since $\psi$, $\varphi$ and $g$ are strictly concave, the first order conditions characterize the points maximizing each function between braces and, with the convention that $(\phi')^{-1}(\infty)=0$ for $\phi=\psi,\varphi,g$, we have :
\begin{equation}
m^*=(g')^{-1}(1/\gamma)~;~
\label{opt2}
e^*=(\varphi')^{-1}\left(\frac{1}{(1-\beta)^+}\right)~;~
a^*=(\psi')^{-1}\left(\frac{1+r_J}{(\delta(1-\beta))^+}\right).
\end{equation}
Let us now consider Situation $2$ and suppose that $(a,\alpha,\beta,\gamma)$ is given. Similarly we obtain that
\begin{equation}
\label{e*m*}
m^*=(g')^{-1}(1/\gamma),~
e^*=(\varphi')^{-1}\left(\frac{1}{(1-\beta)^+}\right).
\end{equation}
When useful to explicit the dependence of these best responses of firm $J$ with respect to the controls of firm $I$, we use the notation $m^*(\gamma)$, $e^*(\beta)$ and $a^*(\beta)$.

To describe the Nash and Stackelberg equilibria, we introduce the continuous mapping $C:\mathbb{R}\to\mathbb{R}$, $B:\mathbb{R}_+\to\mathbb{R}_+$ and $\tilde B:\mathbb{R}\times\mathbb{R}_+\to\mathbb{R}_+$ defined by 

\begin{eqnarray}
 \label{defA} C(\beta)&:=&\frac{1}{v}\ln\mathbb{E}\left[e^{v(1-\beta)(\mu-\varphi( e^*(\beta)))}\right],\\
B(m)&:=&e^{u(Id-b^m)( m)}e^{u(Id-b^e)(e^*(\frac{v}{u+v}))}e^{(u+v) C(\frac{v}{u+v})},\label{defBnash}\\
\label{defB}\tilde B(\beta,m)&:=&e^{u(Id-b^m)(m)}e^{u(Id-b^e)\circ e^*(\beta)}e^{uC(\beta)}\mathbb{E}\left[e^{u\beta(\mu~-\varphi\circ e^*(\beta))}\right].
\end{eqnarray}


\subsection{Nash equilibria}\label{ssec:Nash}
To describe the Nash equilibria, we need the following technical result about the function $g$ :
\begin{lemma}
\label{Lemmalimginfty}
The function $G:m\mapsto \frac{g(m)}{g'(m)}-m$ is continuous, satisfies 
$\lim_{m\to\infty}\frac{g(m)}{g'(m)}-m=+\infty$
and is decreasing from $0$ to $-m_0$ on $[0,m_0]$ (where $m_0$ is defined in (\ref{hypginfty})) and increasing from $-m_0$ to $+\infty$ on $[m_0,+\infty)$ thus admitting an continuous inverse $G^{-1}:[-m_0,+\infty)\to [m_0,+\infty)$.
\end{lemma}



\begin{theorem}
\label{propnash1}
In Situation 1, there exists an infinite number of Nash equilibria, namely the vectors 
$(\hat\alpha, \hat\beta, \hat\gamma, \hat e,\hat m, \hat a)$ satisfying
\begin{eqnarray}
&\hat\beta=\frac{v}{u+v},~\hat e=(\varphi')^{-1}(\frac{u+v}{u}),~\hat a=
(\psi')^{-1}\left(\frac{(u+v)(1+r_J)}{\delta^+ u}\right),~ \hat\gamma=\frac{1}{g'(\hat m)}, 
\label{nash11}\\
&\hat\alpha=C(\hat \beta)+(1+r_J)(D+\hat a)+\hat e-\frac{u}{u+v}\delta\psi(\hat a)-G(\hat m),
\label{nash12}\end{eqnarray}
(where $C$ is defined in (\ref{defA}) 
)
for $\hat{m}$ varying in 
\\
$\hat{\cal M}_1(r_J):=[0,G^{-1}(C(\hat \beta)+(1+r_J)(D+\hat a)+\hat e-\frac{u}{u+v}\delta\psi(\hat a))].$
\\
The corresponding optimal values for firms $J$ and $I$ are respectively $V(0)$ and
$$\hat{I}^1(\hat{m})= b^a(\hat a)-e^{-u\delta\psi(\hat a)} e^{u(1+r_J)(D+\hat a)}B(\hat m),$$
where $B$ is defined in \eqref{defBnash}.
\end{theorem}

 \begin{remark} \label{remmaxnash}
 \begin{itemize}
\item Although there exists an infinite number of Nash equilibria, the controls  $\beta$, $e$ and $a$ are the same in all these equilibria. 
   \item Since $\hat{\mu}-\varphi(\hat e)-\delta\psi(\hat a)\geq 0$, one has  $C(\hat \beta)-\frac{u}{u+v}\delta\psi(\hat a)\geq 0$ so that\\ $[0,G^{-1}((\varphi')^{-1}(\frac{u+v}{u})+D)]\subset\bigcap_{r_J\geq 0}\hat{\cal M}_1(r_J)$.\item It is natural to wonder whether there exists in Situation 1  a Nash equilibrium among the infinite family of such equilibria exhibited in Theorem \ref{propnash1} which maximizes $\hat{I}^1$. This function depends on the Nash equilibrium only through the term $b^m(\hat m)-\hat{m}$ which has to be maximized. The function $\hat m\mapsto b^m(\hat m)-\hat{m}$ being concave, it admits a unique maximum on the interval $\hat{\cal M}_1(r_J)$ where $\hat{m}$ associated with a Nash equilibrium varies. When $[(b^m)']^{-1}(1)\in\hat{\cal M}_1(r_J)$ (which is the case for the numerical example investigated in Section \ref{sec:example}), then $\sup_{\hat{m}\in\hat{\cal M}_1(r_J)}\hat{I}^1(\hat{m})=\hat{I}^1([(b^m)']^{-1}(1))$ and the optimal value $[(b^m)']^{-1}(1)$ of $\hat{m}$ will turn out to be the optimal control $m$ in the Stackelberg equilibria (see Theorems \ref{PropStackelbergpb1} and \ref{PropStackelberg} below).
	\\
	Otherwise, $\sup_{\hat{m}\in\hat{\cal M}_1(r_J)}\hat{I}^1(\hat{m})=\hat{I}^1(G^{-1}(C(\hat \beta)+(1+r_J)(D+\hat a)+\hat e-\frac{u}{u+v}\delta\psi(\hat a)))$.
 \end{itemize}
 \end{remark}

\begin{theorem}
\label{propnash2}
 Let $F(a)=b^a(a)-f((1+r_I)(D+a))$. 
In Situation 2, there exists an infinite number of Nash equilibria namely the vectors $(\hat\alpha, \hat\beta, \hat\gamma, \hat e,\hat m, \hat a)$ defined by
\begin{eqnarray}
\label{NashEqui}
&\hat{m}\geq 0,~\hat\beta=\frac{v}{u+v},~\hat e=(\varphi')^{-1}\Big(\frac{u+v}{u}\Big),~\hat\gamma=\frac{1}{g'(\hat m)},
\\
\label{NashEqui2}
&\hat a\in\arg\max_{a\geq 0} \left[ F(a)-e^{-u\delta\psi(a)}   B(\hat m)\right],
\\\label{NashEqui3}
&\hat\alpha=C(\hat \beta)+\hat e-\frac{u}{u+v}\delta\psi(\hat a)-G(\hat m),
\end{eqnarray}
and such that $C(\hat \beta)+\hat e-\frac{u}{u+v}\delta\psi(\hat a)-G(\hat m)\geq 0$, a condition that is satisfied when $\hat{m}\leq G^{-1}((\varphi')^{-1}(\frac{u+v}{u}))$. 
Moreover, $\hat\alpha+\hat\gamma>0$ and if $\delta\geq 0$, then $\hat a$ is positive and unique for each $\hat m$.
\\
The corresponding optimal values for firms $J$ and $I$ are respectively $V(0)$ and
$$\hat{I^2}(\hat m)=F(\hat a)-e^{-u\delta\psi(\hat a)}
B(\hat m).$$
\end{theorem}

\noindent Let $\hat{\cal M}_2(r_I)$ denote the set of $\hat m\geq 0$ for which there exists $(\hat\alpha, \hat\beta, \hat\gamma, \hat e,\hat a)$ such that $(\hat\alpha, \hat\beta, \hat\gamma, \hat e,\hat m, \hat a)$ is a Nash equilibrium in Situation 2. 

\begin{remark}
\label{remmaxstack}
\begin{itemize}
   \item Notice that the order the different controls are determined is important, since some of them depend on the other ones. Indeed $\hat\beta$ depends on no other control and therefore should be determined first, leading to the value of $\hat e$. Then one should fix $\hat m$, in order to have $\hat\gamma$, which allows then to determine $\hat a$, and once this is done, we can find $\hat\alpha$. Although $\hat\alpha$ and $\hat\gamma$ essentially play the same role, the fact that $\hat\gamma$ only depends on $\hat m$ makes this order important. If one chooses $\hat\alpha$ first, then the determination of $\hat a$ is not clear, since then $\hat a$ depends on $\hat m$, while $\hat m$ depends on $\hat a$ and $\hat\alpha$.
\item We expect that, as in Situation 1, when $[(b^m)']^{-1}(1)\in\hat{\cal M}_2(r_I)$ { (which is the case for the numerical example investigated in Section \ref{sec:example})}, then $\sup_{\hat{m}\in\hat{\cal M}_2(r_I)}\hat{I^2}(\hat{m})=\hat{I}^1([(b^m)']^{-1}(1))$. Indeed, a formal differentiation of $\hat{I^2}(\hat m)$ with respect to $\hat{m}$ leads to 
$(\hat{I^2})'(\hat m)=-e^{-u\delta\psi(\hat a)}
B'(\hat m)$ since, because of the first order optimality condition related to 
\eqref{NashEqui2}, one should not need to take care of the dependence of $\hat a$ on $\hat m$. Moreover, one easily checks that the unique solution to $B'(m)=0$ is $m=[(b^m)']^{-1}(1)$. This is illustrated in Figure \ref{fig:I1(m)andI2(m)}.
\end{itemize}
\end{remark}

We now compare the respective optimal values $\hat I^1(\hat m)$ and $\hat I^2(\hat m)$ for firm $I$ in Situations $1$ and $2$ for the same maintenance effort
 $\hat{m}\in\hat{\cal M}_1(r_J)\cap\left\{\bigcap_{r_I\geq 0}\hat{\cal M}_2(r_I)\right\}$. 
According to Remarks \ref{remmaxnash} and \ref{remmaxstack}, the same value $\hat m=[(b^m)']^{-1}(1)$ is likely to maximize  $ \hat{I}^1(\hat m)$ and    $\hat{I}^2(\hat m)$. Therefore choosing the same maintenance effort $\hat{m}$ for the comparison is natural.
Notice also that, by Theorem \ref{propnash2} and Remark \ref{remmaxnash},  $[0,G^{-1}((\varphi')^{-1}(\frac{u+v}{u}))]\subset\hat{\cal M}_1(r_J)\cap\left\{\bigcap_{r_I\geq 0}\hat{\cal M}_2(r_I)\right\}$. Let $\hat a_1(r_J)=(\psi')^{-1}\left(\frac{(u+v)(1+r_J)}{\delta^+ u}\right)$ (resp. $\hat a_2(r_I)$) denote the value of the control $a$ in the Nash equilibrium with $m=\hat{m}$ in Situation 1 (resp. in Situation 2 when $\delta\geq 0$).

We are going to exhibit cases  in which Situation 1 (meaning outsourcing (respectively Situation 2, meaning debt issuance) is the more profitable for  firm I.

\begin{theorem}
 \label{prop:1stCS}
Let rate $r_J\geq 0$ be fixed and $\hat{m}\in \hat{\cal M}_1(r_J)\cap\left\{\bigcap_{r_I\geq 0}\hat{\cal M}_2(r_I)\right\}$. In case of rate $r_I$ satisfying
  \begin{equation}
 \label{ccl2}
 r_I \leq \frac{ f^{-1}\left[
 {B(\hat m)} e^{-u\delta\psi(\hat a_1(r_J))}\left(e^{u(1+r_J)(D+\hat a_1(r_J))}-1\right)\right]}{D+\hat a_1(r_J)}-1,
 \end{equation}
we have $\hat{I}^2(\hat m)\geq \hat I^1(\hat m)$ and the better contract for firm $I$ is the second one, meaning debt issuance.
 \end{theorem}
Condition (\ref{ccl2}) has a clear economical interpretation. The right-hand side  does not depend on $r_I$. Therefore for a fixed $r_J$, debt issuance is the best  choice  for firm $I$ as soon as its interest rate $r_I$  is small enough. Note the impact  of the function $f$ modeling its debt aversion : the larger $f$, the smaller  the threshold on $r_I$ in  condition (\ref{ccl2}), look at Figure \ref{fig:I1(r_J)I2(r_I)Nash}.

\begin{theorem}
 \label{prop:2dCS}
 We assume $\delta>0.$
Let rate $r_J\geq 0$ be fixed (thus $\hat a_1(r_J)$ is fixed) and $\hat{m}\in \hat{\cal M}_1(r_J)\cap\left\{\bigcap_{r_I\geq 0}\hat{\cal M}_2(r_I)\right\}$. In case of rate $r_I$ satisfying
 \begin{equation}
 \label{ccl1}
(1+r_I)(D+\hat a_2(r_I))\geq f^{-1}\left[
 {B(\hat m)}e^{-u\delta\psi(\hat a_2(r_I))}\left(e^{u(1+r_J)(D+ \hat a_2(r_I))}-1\right)\right] 
 \end{equation}
 and one of the following:
\begin{eqnarray}
\label{CSI^1(a1)(a2)}
& \hspace{-0.5cm} (1+r_I)f'[(1+r_I)(D+\hat a_1(r_J))]&>(b^a)'(\hat a_1(r_J))+u\delta  \psi'(\hat a_1(r_J)) B(\hat m) e^{-u\delta\psi(\hat a_1(r_J))}  ,\nonumber
\\
&\mbox{or}&
\\
&\psi'(\hat a_2(r_I))&>\frac{(u+v)(1+r_J)}{u\delta} \nonumber
\end{eqnarray}
one has $\hat I^1(\hat m)\geq \hat{I}^2(\hat m)$ and 
 the better contract for  firm $I$ is the first one, meaning outsourcing.
 \end{theorem}
The economical interpretation of condition (\ref{CSI^1(a1)(a2)}) is  natural. 
Indeed, the right-hand sides of the inequalities do not depend on $r_I$ whereas the 
left-hand sides are increasing functions of $r_I$. Hence (\ref{CSI^1(a1)(a2)}), leading
 to optimality of outsourcing for firm $I$,  is satisfied as soon as its interest rate
  $r_I$ is large enough. 
Besides, we see that the more convex $f$ is, the smaller is  the threshold on $r_I$ in
 the first inequality of condition (\ref{CSI^1(a1)(a2)}).
 
 Unfortunately, we have not been able to check that the
   condition (\ref{ccl1}) for optimality is satisfied for large $r_I$,
	but Figure \ref{fig:zonesdecisionNash} gives a rule
	of decision between debt issuance and outsourcing.

\subsection{Stackelberg equilibria}
\label{sec:stack}
Depending on Situation 1 or 2 and on the sign of $\delta$, the optimal $\beta$ will be characterized as solution of different equations. To specify those equations, we need to introduce the functions 
\begin{equation}
\label{def:h}
h(\lambda)=\frac{\mathbb{E}\left[\mu e^{\lambda\mu}\right]}{\mathbb{E}\left[e^{\lambda\mu}\right]},
\end{equation}
\begin{equation}\label{def:S}
S(\beta):=\frac{\frac{\beta}{1-\beta}+(b^e)'\circ(\varphi')^{-1}\Big(\frac{1}{(1-\beta)^+}\Big)}{(1-\beta)^2\varphi''\circ(\varphi')^{-1}\Big(\frac{1}{(1-\beta)^+}\Big)}, 
 \end{equation}
{\small 
\begin{equation}
\label{def:Stilde}
\widetilde{S}(\beta):=
\frac{1+r_J}{\delta(1-\beta)^2 (\psi'')((\psi')^{-1}\Big(\frac{1+r_J}{(\delta(1-\beta))^+}\Big)) } \left(  (1+r_J)\frac{\beta}{1-\beta} +  \frac{(b^a)'\circ(\psi')^{-1}\Big(\frac{1+r_J}{(\delta(1-\beta))^+}\Big)}{ u  e^{ (Id-b^m) ) ((b^m)')^{-1}(1)} k(\beta)  }      \right),
\end{equation}}
where $k(\beta)$ is a positive function of $\beta$ defined  as follows,
{\small 
\begin{equation}
\label{defk}
k(\beta):=e^{u(Id-b^e)\circ e^*(\beta)}   e^{-u \delta \psi \circ a^*(\beta)}e^{u(1+r_J)(D+a^*(\beta))}   e^{uC(\beta)}\mathbb{E}\left[e^{u\beta\big(\mu-\varphi\circ e^*(\beta)\big)}\right]
\end{equation}}
with $e^*(\beta)$ and $a^*(\beta)$ defined in (\ref{opt2}) and $C(\beta)$ in (\ref{defA}).\\

We consider the  following equations
\begin{equation}
\label{eqbetaStackelberg}
h(u\beta)-h(v(1-\beta))={S} (\beta),
\end{equation}

\begin{equation}
\label{eqbetaStackelbergpb1positif}
h(u\beta)-h(v(1-\beta))=S(\beta)+\widetilde{S} (\beta),
\end{equation}

\begin{equation}
\label{eqbetaStackelbergpb1negatif}
h(u\beta)-h(v(1-\beta))=\widetilde{S} (\beta).
\end{equation}

\begin{theorem}
\label{PropStackelbergpb1}
In Situation 1, there exists at least one Stackelberg equilibrium with firm $I$ as the leader. Moreover, if there exists   a Stackelberg equilibrium 
$(\hat e, \hat m, \hat a, \hat\alpha, \hat\beta,\hat\gamma)$
with $\hat\alpha>0$, then it is characterized by :
 $$\hat e=(\varphi')^{-1}\Big(\frac{1}{(1-\hat\beta)^+}\Big),~
 \hat a=(\psi')^{-1}\left(\frac{1+r_J}{(\delta(1-\hat\beta))^+}\right),
 \hat m=\big[(b^m)'\big]^{-1}(1),~\hat\gamma=1/g'(\hat m),$$
  $$\hat\alpha= C(\hat\beta)+(1+r_J)(D+\hat a)+\hat e-(1-\hat\beta)\delta\psi(\hat a)-G(\hat m).$$
   If $\delta>0$ then 
$\hat \beta$ is a solution of (\ref{eqbetaStackelbergpb1positif}) and is less than $\frac{v}{u+v}$. 
\\
If $\delta<0$ then either $\hat \beta$ is less than $\frac{v}{u+v}$ and solves (\ref{eqbetaStackelberg}) or $\hat \beta$ is larger than one and solves (\ref{eqbetaStackelbergpb1negatif}).\\
The corresponding optimal values for firms $J$ and $I$ are respectively $V(0)$ and $I^1(\hat\alpha,\hat\beta,\hat\gamma)$.
\end{theorem}

\begin{theorem}
\label{PropStackelberg}
In Situation 2,  there exists at least one Stackelberg equilibrium with firm $I$ as the leader.
Moreover, if there exists   a Stackelberg equilibrium 
$(\hat e, \hat m, \hat a, \hat\alpha, \hat\beta,\hat\gamma)$
with $\hat\alpha>0$, then it satisfies:
\\
$\hat\beta$ is a solution of (\ref{eqbetaStackelberg}), $\hat e=(\varphi')^{-1}\Big(\frac{1}{(1-\hat\beta)^+}\Big)$, $\hat m=\big[(b^m)'\big]^{-1}(1)$, $\hat\gamma=1/g'(\hat m)$, 
$\hat a\in\arg\max_{a\geq 0} F(a)-e^{-u\delta\psi(a)}\tilde B(\hat\beta,\hat m)$ and
 $\hat\alpha=C(\hat\beta)+\hat e-(1-\hat\beta)\delta\psi(\hat a)-G(\hat m)$,
 where the mappings $C$ and $\tilde B$ are defined by (\ref{defA})-(\ref{defB}).
In particular, $\hat\beta<\frac{v}{u+v}$. 
\\
Moreover, if $\delta\geq 0$, then $\arg\max_{a\geq 0} F(a)-e^{-u\delta\psi(a)}\tilde B(\hat\beta,\hat m)$ is a singleton and $\hat a>0$.\\
The corresponding optimal values for firms $J$ and $I$ are respectively $V(0)$ and $I^2(\hat a,\hat\alpha,\hat\beta,\hat\gamma)$.\end{theorem}

Proposition \ref{eq5.5} below states that Equation (\ref{eqbetaStackelberg}) which appears in the characterization of $\hat{\beta}$ when $\hat \alpha>0$ in both Situations 1 and 2 always admits a solution.

 An analytical comparison is not so easy, but Figures 4 and
 \ref{fig:zonesdecisionStackelberg} allow a numerical comparison between debt issuance and outsourcing.

\subsection{Incomplete information}
\label{incompl}
In this section we consider the previous equilibrium problems when the firms do not have a perfect knowledge of the preferences of the other firm. More precisely, we still assume that the firms' utility functions are $U(x)=-e^{-ux}$ and $V(x)=-e^{-vx}$ respectively, but firm $I$ perceives $v$ as 
 a $(0,+\infty)$-valued random variable  with known distribution and independent from $\mu$ that we denote ${\cal V}$ and firm $J$ perceives $u$ as a random variable with known distribution and independent from $\mu$
that we denote ${\cal U}$. According to Section (\ref{bestrespj}),  firm $J$ optimal controls are functions of the 
controls $\beta,\gamma$ fixed by firm $I$ that do not depend on the risk aversion parameters
 $u$. Therefore, equations (\ref{opt2}) and (\ref{e*m*})    still hold in incomplete information
  and  incomplete information on the risk aversion parameter $u$ has no impact on the equilibria. 
In contrast, the uncertainty on the parameter $v$ has an impact as the acceptation of the
 contract by firm $J$ depends on it.  To model the social need of the investment, we 
 introduce a (social) penalty $p\in\mathbb{R}\cup\{+\infty\}$ that firm $I$ gets if firm 
 $J$ does not accept the contract.

\subsubsection{Stackelberg equilibrium, firm $I$ is leader}\label{stackincominf}
We first introduce
the events ${\cal A}^i,~i=1,2:$  ``firm $J$ accepts the contract'' in Situation $i$.
\\
The optimization problem for firm $I$ is
\begin{equation}
\label{eqinfoincomplete1}
 u_I^1:=-p\vee\sup_{(\alpha,\beta,\gamma)}\{\tilde I^1(\alpha,\beta,\gamma)\pr\left({\cal A}^1(\alpha,\beta,\gamma)\right)-p(1-\pr\left({\cal A}^1(\alpha,\beta,\gamma)\right))\},
 \end{equation}
in Situation 1 and in Situation 2, it becomes :
\begin{equation}
\label{eqinfoincomplete2}
 u_I^2:=-p\vee\sup_{(a,\alpha,\beta,\gamma)}\{\tilde I^2(a,\alpha,\beta,\gamma)\pr\left({\cal A}^2(a,\alpha,\beta,\gamma)\right)-p(1-\pr\left({\cal A}^2(a,\alpha,\beta,\gamma)\right))\}.
 \end{equation}
 The functions 
\begin{eqnarray*}
  \hspace*{-0.5cm} 
  \tilde I^1(\alpha,\beta,\gamma)&:=&
b^a (a^*(\beta))
- \mathbb{E}\left[ 
e^{-u\big([b^m-\gamma g](m^*(\gamma))+[b^e+\beta\varphi](e^*(\beta))
-\alpha-\beta(\mu-\delta \psi(a^*(\beta)))\big)}\right ],
\\
\tilde I^2(a,\alpha,\beta,\gamma)&:=&F(a)-
\mathbb{E}\left[e^{-u\big([b^m-\gamma g](m^*(\gamma))+[b^e+\beta\varphi](e^*(\beta))-\alpha-\beta(\mu-\delta \psi(a))\big)}\right],
\end{eqnarray*}
where $e^*,~m^*$ and $a^*$ have been defined in (\ref{opt2}),
are the social gain that firm $I$ respectively gets in Situations $1$ and $2$ if firm $J$ accepts the contract. Notice that the supremum is taken with $-p$ to model the possibility for firm $J$ not to enter the game and that $p=+\infty$ corresponds to the case where firm $I$ absolutely wants that firm $J$ accepts the contract.
\\
In order to characterize the acceptance set ${\cal A}^i$, we introduce
\begin{align}
\tilde J^1(v,\alpha,\beta,\gamma)&:=-e^{v(1+r_J)(D+ a^*(\beta))}\notag\\&\times\E e^{-v\big(\alpha+(\beta-1)(\mu-\delta\psi(a^*(\beta))-\varphi(e^*(\beta)))-e^*(\beta)+[\gamma g-Id](m^*(\gamma)) \big)}  
\label{deftJv1}
 \end{align}
 and
 \begin{equation}
 \label{deftJv2}
\tilde J^2(v,a,\alpha,\beta,\gamma):=
-\E e^{-v\big(\alpha
+(\beta-1)(\mu-\delta\psi(a)-\varphi(e^*(\beta)))-e^*(\beta)+
[\gamma g-Id](m^*(\gamma))\big)}.
\end{equation}
Firm $J$ accepts the contract if and only if $-{\tilde J^i}({\cal V},.)\leq 1$, thus ${\cal A}^i(.)=\{-{\tilde J^i}({\cal V},.)\leq 1\}$.\\
We define  the value function of the problem with complete information that firm $J$'s risk aversion is equal to $v$ $$u^1(v):=\sup_{\{(\alpha,\beta,\gamma)\in\mathbb{R}_+\times\mathbb{R}\times\mathbb{R}_+:- \tilde J^1(v,\alpha,\beta,\gamma)\leq 1\}}\tilde I^1(\alpha,\beta,\gamma)$$
$$u^2(v):=\sup_{\{(a,\alpha,\beta,\gamma)\in\mathbb{R}_+\times\mathbb{R}_+\times\mathbb{R}\times\mathbb{R}_+:-\tilde J^2(v,a,\alpha,\beta,\gamma)\leq 1\}}\tilde I^2(a,\alpha,\beta,\gamma)$$
These value functions are  respectively obtained for the Stackelbreg equilibria given in Theorems \ref{PropStackelbergpb1} and \ref{PropStackelberg}. We have the following result:
\begin{theorem}
\label{Propinfoincomplete}
Let \begin{equation}
\label{Pbequivalent}
 w^i_I:=-p\vee \sup_{v>0}\{u^i(v)\pr({\cal V}\leq v)-p\big(1-\pr({\cal V}\leq v)\big)\}.
 \end{equation}
We have $w^i_I\leq u^i_I$ and when either $p<+\infty$ or $\exists v\in (0,+\infty),\;\pr({\cal V}>v)=0$ then $w^i_I=u^i_I$.
\end{theorem}
Theorem \ref{Propinfoincomplete} has an important interpretation. Indeed, it means
 that in order to solve (\ref{eqinfoincomplete1}) or (\ref{eqinfoincomplete2}), firm $I$
  first solves its problem for any given $v$ as if the information was complete or in 
  other words as in Section \ref{sec:stack}, and then "chooses" the level $v$ that would 
  bring the greatest social expectation in (\ref{Pbequivalent}).

\vspace{2mm}

\begin{theorem}
\label{Th:stackincomplet} 
Let $v_0:=\inf\{v>0:\pr({\cal V}\leq v)>0\}$. If $\lim_{v\to v_0^+}u^i(v)\leq -p$ then the fact 
that the two firms do not enter into any contract is a Stackelberg equilibrium in Situation~$i$.
\\
Otherwise, if $v_1:=\sup\{v>0:\pr({\cal V}>v)>0\}<+\infty$ then the optimization problem
 (\ref{Pbequivalent}) has a solution $v^\star\in (0,v_1]\cap [v_0,v_1]$ (equal to $v_1$ when $p=+\infty$) and any Stackelberg equilibrium for the problem with complete information and risk aversion $v^\star$ for firm $J$ is a Stackelberg equilibrium for the problem with incomplete information.
 \end{theorem}

\subsubsection{Nash equilibrium}
We did not succeed in finding sufficient conditions for the existence of  a Nash equilibrium with incomplete information. Nevertheless, we obtain necessary conditions that are similar for both situations:

\begin{theorem}
\label{prop:Nashincomplet}
Assume existence of a Nash equilibrium $\hat{c}=(\hat{a},\hat{\alpha},\hat{\beta},\hat{\gamma},\hat{e},\hat{m})$ such that the value for firm $I$ is greater than $-p$ and let
$\hat v:=\sup\{v>0:-J(v,\hat{c})\leq 1\}$
with $J(v,\hat{c}) $  defined (using $C^{op}=\mu-\varphi(e)-\delta\psi(a)$) in Situations 1 and  2 respectively as
$$J(v,a,\alpha,\beta,\gamma,e,m):=-\mathbb{E}\left[ e^{-v\big(\alpha+(\beta-1)C^{op}-e+\gamma{g}(m)-m -(1+r_J)(D+ a)\big)}\right],$$
$$J(v,a,\alpha,\beta,\gamma,e,m):=-\mathbb{E}\left[e^{-v\big(\alpha+(\beta-1)C^{op}-e+\gamma{g}(m)-m\big)}\right].$$
 Then $\hat v>0$,
$\hat e=(\varphi')^{-1}\left(\frac{1}{(1-\hat \beta)^+}\right),~\hat \gamma = \frac{1}{g'(\hat m)}$ and in Situation 1, 
$\hat a=(\psi')^{-1}\left(\frac{1+r_J}{(\delta(1-\hat\beta))^+ }\right)$. 
\\
If $\hat v<+\infty$, then $\hat{c}$ is a Nash equilibrium for the problem with complete information and risk aversion $\hat v$ for firm $J$ and for each $v<\hat v$, $\pr({\cal V}\leq v)<\pr({\cal V}\leq \hat v)$. 
\\
If $\hat v=+\infty$, then for each $v\in (0,+\infty)$, $\pr({\cal V}\leq v)<1$.
\end{theorem}
\begin{remark}
   If there is a vector $(v_1,\cdots,v_n)$ of elements of $(0,+\infty)$ such that 
   \\
   $\sum_{k=1}^n\pr({\cal V}=v_k)=1$, one deduces that if there exists a  Nash equilibrium for the problem with incomplete information, then $\exists i$ such that $\hat v=v_i$.
\end{remark}

\section{Numerical Example}
\label{sec:example}
We investigate a numerical example to better quantify and compare the two different situations (Situation 1: outsourcing, Situation 2 : debt issuance) and the
 two equilibria (Nash and Stackelberg equilibria). We have chosen the following numerical values: 
\begin{itemize}
\item the risk aversion parameters are $u = v = 1$, 
\item  the impact of the quality of the investment on the operational costs is $\delta = 1$,
\item the minimal investment is $D = 1$,
\item the random cost $\mu$  follows  a uniform distribution on the interval $[\underline{\mu} ,\overline{\mu}] $, with $\underline{\mu} = 1$ and $\overline{\mu} = 2$,
\item the aversion to debt issuance  is $f(x) = e^{1.6 x} - 1$  (or $e^{1.7 x} - 1$),
\item the benefits functions are $b^a(x) = \frac{x}{(1+x)}$ and $b^m(x) = b^e(x) = \sqrt{x}$,
\item the impact of the maintenance cost $m$ on the rent $R$  is modeled by the function $g(x) = \sqrt{x}$, and the 
 impact of the efforts $a$ and $e$ on the operational costs   by the functions 
$\phi(x) = \psi(x) = \frac{\underline{\mu}}{3(1+\pi/4)} (  \sqrt{x} {\bf 1}_{x \leq 1} + (\arctan(x) + 1 - \pi/4) {\bf 1}_{x > 1}      )$ 
(thus satisfying the condition $\varphi(\infty)+\delta^+\psi(\infty)< \underline{\mu}$),
\item the interest rates $r_I$ and  $r_J$ over the period $[0,T]$ (when borrowing $1$ initially, firm $K\in\{I,J\}$ has to reimburse $(1+r_K)$ on the time-interval $[0,T]$) vary in the interval  $[0,1[$.
\end{itemize}

\subsection{Nash equilibrium}
\label{fig:Nash}
\begin{figure}[ht]
    \centering
    \includegraphics[width=0.8\textwidth]{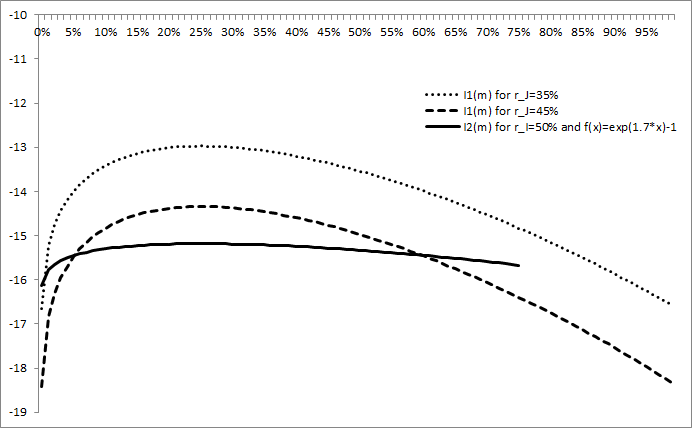}
    \caption{Value functions $\hat{I}^1$ and $\hat{I}^2$ as functions of $\hat{m}$}
    \label{fig:I1(m)andI2(m)}
\end{figure}

\paragraph{Dependency on the  maintenance costs $\hat{m}$ in Nash equilibrium}
We now investigate the Nash equilibria in both situations. As stated in Theorem \ref{propnash1} and Theorem \ref{propnash2}, the optimal value functions $\hat{I}^1$ and $\hat{I}^2$ depend on the optimal maintenance costs $\hat{m}$. 
 Figure \ref{fig:I1(m)andI2(m)} gives an insight of this dependency, for two different values of $r_J$ in Situation 1 and for a larger value of $r_I$ in Situation 2.
We notice that outsourcing is worthless for average maintenance costs (between $5 \% $ and $50 \% $ of $D$ for $r_J=45\%$) 
We also observe that the smaller $r_J$, the larger is the interval of values $\hat{m}$ for which outsourcing is better than debt issuance.

The  maintenance costs $\hat{m}$ numerically optimizing the value functions $\hat{I}^1$ and $\hat{I}^2$  is $[(b^m)']^{-1}(1)=\frac{1}{4}$ (which is also the maintenance costs in Stackelberg equilibrium, whatever the situation) as expected from Remarks \ref{remmaxnash} and \ref{remmaxstack}. Therefore, in the forthcoming figures, the Nash equilibrium is computed for this optimal maintenance costs $\hat{m}=\frac{1}{4}$.
\newpage

\begin{figure}[ht]
    \centering
    \includegraphics[width=0.8\textwidth]{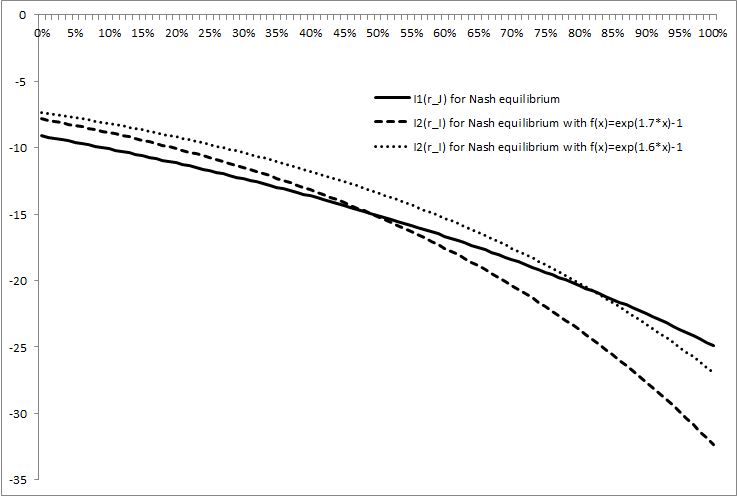}
    \caption{$\hat{I}^1(r_J)$ and $\hat{I}^2(r_I)$ in Nash equilibrium}
    \label{fig:I1(r_J)I2(r_I)Nash}
\end{figure}

\paragraph{Dependency on the interest rates $r_I$ and $r_J$ in Nash equilibrium}
Figure \ref{fig:I1(r_J)I2(r_I)Nash} gives, for Nash equilibrium,  the optimal value function $\hat{I}^1$ in Situation 1 (outsourcing) as a function of $r_J \in [0,1[$,  
and $\hat{I}^2$ in Situation 2 (debt issuance) as a function of $r_I \in [0,1[$ and for two differents functions of debt aversion ($f(x) = e^{1.6 x} - 1$  or  $e^{1.7 x} - 1$).
We notice that the greater the debt aversion function $f$, the smaller the value of $r_I=r_J$ at which it becomes more favourable to outsource ($50\%$ for $f(x) = e^{1.7 x} - 1$,   $80\%$  for $e^{1.6 x} - 1$).

\begin{figure}[ht]
    \centering
    \includegraphics[width=0.5\textwidth]{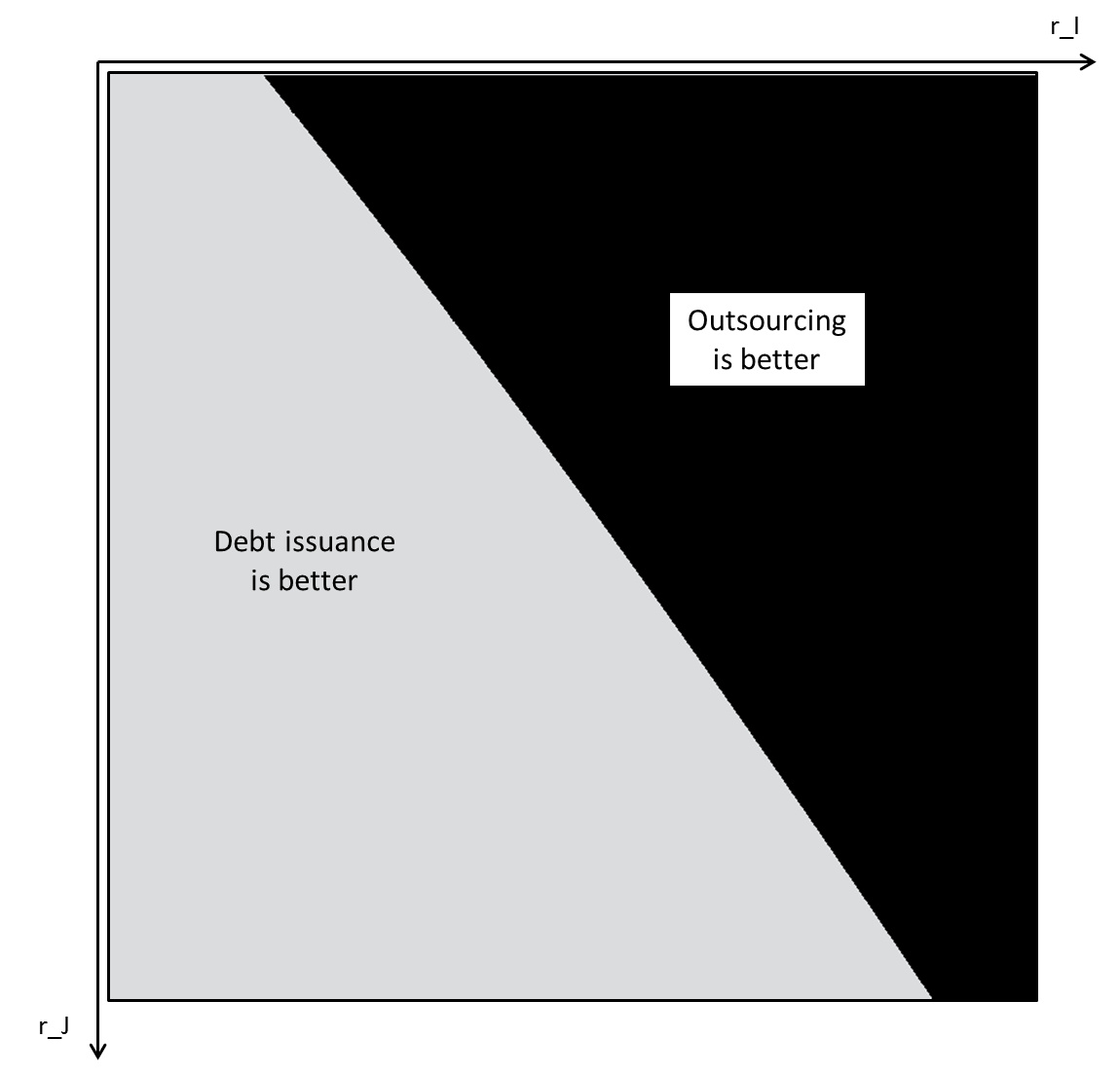}
    \caption{Decision areas in Nash equilibrium}
    \label{fig:zonesdecisionNash}
\end{figure}

\paragraph{Outsourcing or not in Nash equilibrium?}
Figure \ref{fig:zonesdecisionNash} gives a decision criterion of outsourcing or not, function of $r_I\in[0,1[$ in the $x$-axis and $r_J\in[0,1[$ in the $y$-axis. The grey area corresponds to the region where it is optimal to issue debt, while the black area corresponds to the region where it is optimal to outsource. As expected, if it is optimal to outsource for a given couple $(r_I,r_J)$, then it remains optimal to outsource for all couples $(r_I',r_J)$ with $r_I'>r_I$. If it is optimal to issue debt for a given couple $(r_I,r_J)$, then it remains optimal to issue debt for all couples $(r_I,r'_J)$ with $r_J'>r_J$.

\newpage

\subsection{Stackelberg equilibrium}
\label{figStackelberg}


\paragraph{Dependency on the interest rates $r_I$ and $r_J$ in Stackelberg equilibrium}
Figure 4 below gives, for Stackelberg equilibrium,  the optimal value function $\hat{I}^1$ in Situation 1 (outsourcing) as a function of $r_J \in [0,1[$,  and $\hat{I}^2$ in Situation 2 (debt issuance) as a function of $r_I \in [0,1[$ and for two different functions of debt aversion ($f(x) = e^{1.6 x} - 1$  or  $e^{1.7 x} - 1$).
\begin{figure}[ht]
    \begin{center}
    \includegraphics[width=0.8\textwidth]{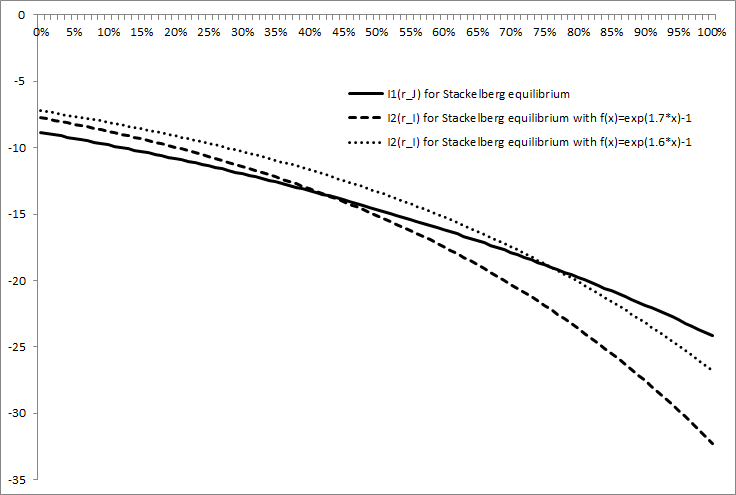}
  \caption{$\hat{I}^1(r_J)$ and $\hat{I}^2(r_I)$ in Stackelberg equilibrium}
\end{center}
\label{fig:**}
\end{figure}
The conclusions are the same as the ones in Figure \ref{fig:I1(r_J)I2(r_I)Nash} for Nash equilibrium, and
we notice that the value functions $\hat{I}^1$ and $\hat{I}^2$ are  slightly greater in Stackelberg equilibrium than  in Nash equilibrium. Moreover, the range of interest rate $r$ for which outsourcing is  more favourable is slightly wider in Stackelberg equilibrium than in Nash equilibrium ($43\%$ for $f(x) = e^{1.7 x} - 1$,   $73\%$  for $e^{1.6 x} - 1$). Indeed, being a leader in Stackelberg equilibrium, firm $I$ has a more favourable position to get better conditions to outsource, in comparison with Nash  equilibrium.

\begin{figure}[ht]
    \centering
    \includegraphics[width=0.5\textwidth]{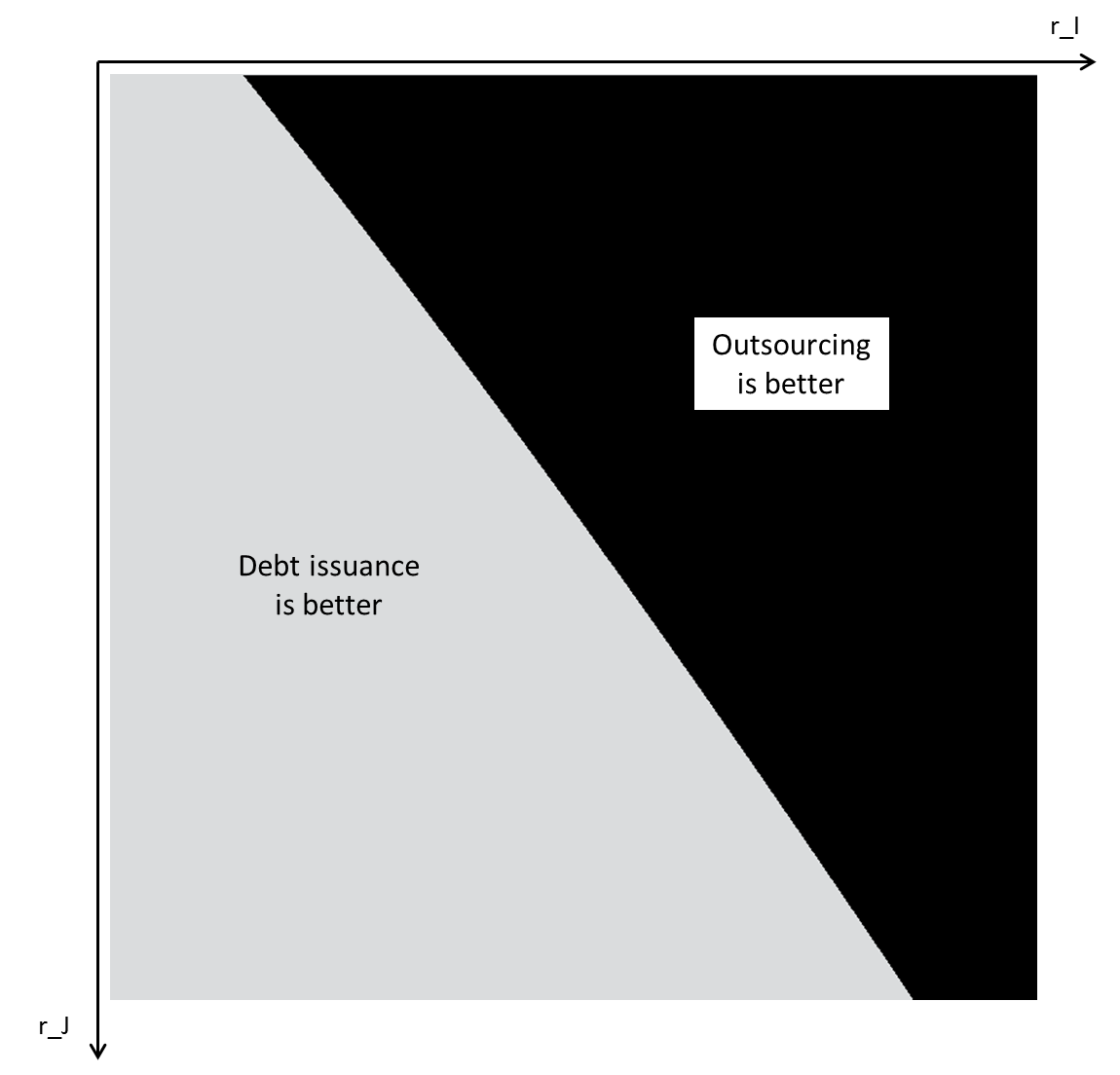}
    \caption{Decision areas in Stackelberg equilibrium}
    \label{fig:zonesdecisionStackelberg}
\end{figure}

\paragraph{Outsourcing or not in Stackelberg equilibrium?}
Figure \ref{fig:zonesdecisionStackelberg} gives a decision criteria of outsourcing or not, function of $r_I$ in the $x$-axis and $r_J$ in the $y$-axis. The grey area corresponds to the region where it is optimal to issue debt, while the black area corresponds to the region where it is optimal to outsource. 
In comparison with Figure \ref{fig:zonesdecisionNash} dedicated to Nash equilibria, the boundary between the outsourcing and debt issuance regions has the same shape but the outsourcing region is slighly bigger.

\newpage

\subsection{Comparison Nash/Stackelberg equilibria}
\begin{figure}[ht]
    \centering
    \includegraphics[width=0.8\textwidth]{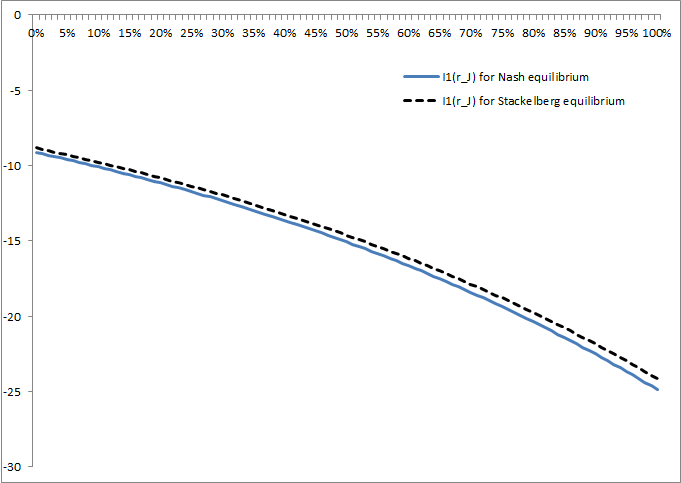}
    \caption{Value function $\hat{I}^1$ (outsourcing) in Nash and Stackelberg equilibria}
    \label{fig:I1NashvsStackelberg}
\end{figure}

\begin{figure}[ht]
    \centering
    \includegraphics[width=0.8\textwidth]{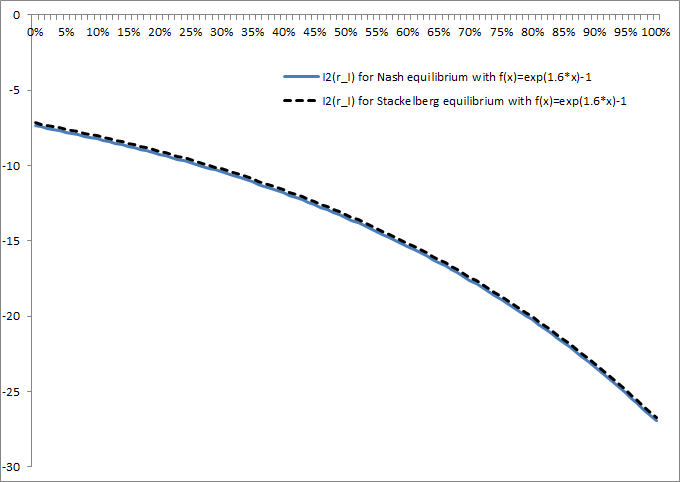}
    \caption{Value function $\hat{I}^2$ (debt issuance) in Nash and Stackelberg equilibria}
    \label{fig:I2NashvsStackelberg}
\end{figure}

We first compare  the value functions  of firm $I$ between Nash and Stackelberg equilibria, in the situation of outsourcing the debt (Figure \ref{fig:I1NashvsStackelberg})
and in the situation of debt issuance (Figure \ref{fig:I2NashvsStackelberg}). In both situations, the value function is higher in Stackelberg equilibrium than in Nash equilibrium, which can be interpreted by the fact that  firm $I$ has a leader position in Stackelberg equilibrium. The difference is a little less significant in Situation 2 than in  Situation 1.

\newpage

\section{Proofs in Situation 1}
\label{sec:proof1}

\subsection{Best responses in Situation 1}
\label{sit1}


Let $a$, $e$ and $m$ be given and constant. Then we get the following optimization problem for firm $I$:
$$
\sup_{\alpha,\beta,\gamma}\mathbb{E}\left[b^a(a)
-  e^{-u \big(b^m(m)+b^e(e)-\alpha-\beta(\mu-\varphi(e)-\delta \psi(a))-\gamma g(m) \big) }\right]
$$
such that $\alpha\geq 0,~\gamma\geq 0,$ and using $V(0)=-1,$
%
\begin{equation}
\label{constr}
\mathbb{E}\left[ e^{-v  \big(\alpha+(\beta-1)(\mu-\varphi(e)-\delta\psi(a))-e+\gamma g(m)-m -(1+r_J)(D+a)\big) }\right]\leq 1.
\end{equation}

\begin{proposition}
\label{PropOptimumFirm1bis}
Let $a\geq 0$,  $e\geq 0$ and $m\geq 0$ be given and constant. Then there exist optimal triplets $(\alpha,\beta,\gamma)$ for the above problem. Moreover $(\alpha,\beta,\gamma)$ is optimal if and only if it satisfies: 
$\beta=\beta^*:=\frac{v}{u+v}\in(0,1)$ 
\\
and 
$\alpha+\gamma {g}(m)=C_e(\beta^*)-\frac{u}{u+v}\delta\psi(a) +e+m+(1+r_J)(D+a)$
 with $\alpha\geq 0$ and $\gamma\geq 0$, where
\begin{equation}
\label{defAbis} C_e(\beta):=\frac{1}{v}\ln\mathbb{E}\left[e^{v(1-\beta)(\mu-\varphi(e))}\right].
 \end{equation}
\end{proposition}
\proof
We first need the following lemmas :
\begin{lemma}
\label{lem:h}
The function $h(\lambda)=
\frac{\mathbb{E}\left[\mu e^{\lambda\mu}\right]}{\mathbb{E}\left[e^{\lambda\mu}\right]}$ 
is increasing, thus the equation $h(u\beta)=h(v(1-\beta))$ admits the unique solution $\beta^*:=\frac{v}{u+v}\in(0,1)$.
\end{lemma}
\proof

We  compute using Cauchy-Schwarz inequality
\begin{equation*} \Big(\mathbb{E}e^{\lambda\mu}\Big)^2h'(\lambda)=\mathbb{E}\mu^2e^{\lambda\mu}\mathbb{E}e^{\lambda\mu}-\Big(\mathbb{E}\mu e^{\lambda\mu}\Big)^2>0,
\end{equation*}
The inequality is strict since $\mu$ is not constant $d\pr$ a.s.. Therefore $\beta$ satisfies $h(u\beta)=h(v(1-\beta))$  if and only if $u\beta=-v(\beta-1)$, so that the only solution is $\beta^*:=\frac{v}{u+v}\in(0,1)$.\\
\ep

\begin{lemma}
   The functions $C_e(\beta)$ and $C(\beta)=C_{e^*(\beta)}(\beta)$ defined in \eqref{defAbis} and \eqref{defA} are such that
\begin{align}
   &\forall \beta<1,\;\forall a,e\geq 0,\;\forall\delta\in{\mathbb R},\;C_e(\beta)-\delta^+(1-\beta)\psi(a)>0\label{minoce}\\
&\forall \beta<1,\;\forall a\geq 0,\;\forall\delta\in{\mathbb R},\;C(\beta)-\delta^+(1-\beta)\psi(a)>0\label{minoc}\\
&\forall \beta\in{\mathbb R},\;
\label{minC}
C(\beta)\geq (1-\beta)\Big(  \mathbb{E}[\mu]~  -\varphi(e^*(\beta))\Big)\end{align}
\end{lemma}
\proof
Since $v>0$, $C_e(\beta)-\delta^+(1-\beta)\psi(a)=\frac{1}{v}\ln\E\left[e^{v(1-\beta)(\mu-\varphi(e)-\delta^+\psi(a))}\right]$ and $\mu-\varphi(e)-\delta^+\psi(a)>0$ a.s., the two first statements are consequences of the monotonicity of $x\mapsto\exp v(1-\beta)x$. The last one is a consequence of Jensen's inequality.
\ep
\\

Let $a$, $e$ and $m$  be given and $C^{op}=\mu-\varphi(e)-\delta\psi(a),$ we introduce:
\begin{eqnarray*}
&K(\alpha,\beta,\gamma):=-\mathbb{E}\left[e^{-u\big(b^m(m)+b^e(e)-\alpha-\beta C^{op}-\gamma{g}(m)\big)}\right]
\\
&E:=\Big\{(\alpha,\beta,\gamma)\in\mathbb{R}^3;~\mathbb{E}\left[e^{-v\big(\alpha+(\beta-1)C^{op}-e+\gamma{g}(m)-m -(1+r_J)(D+a)\big)}\right]\leq 1\Big\}.
\end{eqnarray*}
We will first solve the problem of maximization of $K$ on $E$, forgetting about the constraints $\alpha\geq 0$ and $\gamma\geq 0$, and we will then see that it allows to solve the original constrained problem. Let us therefore consider the following problem: $\sup_{(\alpha,\beta,\gamma)\in E}K(\alpha,\beta,\gamma)$. Since $K$ is concave and $E$ is a closed convex set, the first order conditions for the Lagrangian associated to this problem are also sufficient conditions. The Lagrangian is given by:
\begin{eqnarray*}
&&L(\alpha,\beta,\gamma,\lambda):=-\mathbb{E}\left[e^{-u\big(b^m(m)+b^e(e)-\alpha-\beta C^{op}-\gamma{g}(m)\big)}\right]
\\
&&-\lambda \Big(\mathbb{E}\left[
e^{-v\big(\alpha+(\beta-1)C^{op}-e+\gamma{g}(m)-m -(1+r_J)(D+a)\big)}\right]-1\Big).
\end{eqnarray*}
Hypothesis (\ref{hypmu}) implies that $L$ is differentiable and  
the following are null:
{\small{
\begin{eqnarray}
\label{Eqalphabis}
\frac{\partial L}{\partial \alpha}&=&-\mathbb{E}\left[u e^{-u()}\right]+
\lambda\mathbb{E}\left[ve^{-v()}\right],
\\
\label{Eqbetabis}
\frac{\partial L}{\partial \beta}&=&-\mathbb{E}\left[u e^{-u()}C^{op}\right]+
\lambda\mathbb{E}\left[ve^{-v()}C^{op}\right],
\\
\label{Eqgammabis}
\frac{\partial L}{\partial \gamma}&=&-\mathbb{E}\left[u e^{-u()}g(m)\right]+
\lambda\mathbb{E}\left[ve^{-v()}g(m)\right].
\end{eqnarray}
}}
\noindent Since $g(m)$ is a constant, equation (\ref{Eqalphabis}) implies equation (\ref{Eqgammabis}).
Furthermore, since
 $C^{op}=\mu-\varphi(e)-\delta\psi(a)>0$,  then $\lambda>0$ and the constraint is always binding. This is natural since $K$ is decreasing with respect to $\alpha$, $\beta$ and $\gamma$, while 
 $$(\alpha,\beta,\gamma)\mapsto -
 \mathbb{E}\left[e^{-v\big(\alpha+(\beta-1)C^{op}-e+\gamma g(m)-m-(1+r_J)(D+a)\big)}\right]$$
  is increasing with respect to $\alpha$, $\beta$ and $\gamma$.
   Therefore, at an interior point of $E$ denoted $(\alpha,\beta,\gamma)$, for sufficiently small $\eps>0$, for example $(\alpha,\beta-\eps,\gamma)$ is still in $E$, while 
   $K(\alpha,\beta-\eps,\gamma)>K(\alpha,\beta,\gamma)$, so that $(\alpha,\beta,\gamma)$ cannot be a maximum of $K$.\\
Therefore we also have:
 \begin{equation}
\label{Eqcontraintebis}
\mathbb{E}\left[e^{-v\big(\alpha+(\beta-1)C^{op}-e+\gamma{g}(m)-m-(1+r_J)(D+a)\big)}\right]=1.
\end{equation}
Combining (\ref{Eqalphabis}) and (\ref{Eqbetabis}), we get:
$$
 \mathbb{E}\left[\mu e^{u\beta\mu}\right]
  \mathbb{E}\left[e^{-v(\beta-1)\mu}\right]=
\mathbb{E}\left[e^{u\beta\mu}\right]
\mathbb{E}\left[\mu e^{-v(\beta-1)\mu}\right].
$$
This equation is equivalent to $h(u\beta)-h(v(1-\beta)=0$ which admits the unique solution
$\beta^*=\frac{v}{u+v}\in(0,1)$ (cf. Lemma \ref{lem:h}).\\We have
\begin{equation}
   \mathbb{E}\left[e^{v(1-\beta^*)(\mu~-\varphi(e)-\delta\psi(a))}\right]=e^{v(C_e(\beta^*)-\frac{u}{u+v}\delta\psi(a))}\label{newexpr},
\end{equation} 
which together (\ref{Eqcontraintebis}) yields the following necessary and sufficient condition for optimality: 
 \begin{equation}
 \label{NSCopt}
 \alpha+\gamma g(m)=C_e(\beta^*)-\frac{u}{u+v}\delta\psi(a)+e+m+(1+r_J)(D+a) .
 \end{equation}
Since $1-\beta^*>0$, by \eqref{minoce}, $C_e(\beta^*)-\frac{u}{u+v}\delta\psi(a)>0$, $C_e(\beta^*)-\frac{u}{u+v}\delta\psi(a)+e+m+(1+r_J)(D+a)>0$ and the set
\\
 $\{(\alpha,\gamma)\in[0,+\infty)^2,~\alpha+\gamma g(m)=C_e(\beta^*)-\frac{u}{u+v}\delta\psi(a)+e+m+(1+r_J)(D+a)\}$ is not empty. Therefore the optimal $(\alpha,\beta,\gamma)$ for the problem:
$$\sup_{(\alpha,\beta,\gamma)\in E\cap(\mathbb{R}_+\times\mathbb{R}\times\mathbb{R}_+)}K(\alpha,\beta,\gamma),$$
are exactly the elements of 
\begin{eqnarray*} \big\{(\alpha,\beta,\gamma)\in\mathbb{R}_+\times\mathbb{R}\times\mathbb{R}_+;~\beta&=&\frac{v}{u+v},\\
\alpha+\gamma g(m)&=&C_e(\beta)-\frac{u}{u+v}\delta\psi(a)+e+m +(1+r_J)(D+a)\big\}.
\end{eqnarray*}
\ep

\subsection{Nash equilibrium in Situation 1}

\noindent{\bf Proof of Lemma \ref{propnash1}}
The function $g$ is assumed to be increasing, strictly concave and such that $g'(0)=\infty$, $g'(\infty)=0$ and (\ref{hypginfty}) holds.
Since $g$ is increasing and concave, we compute for any $x\in[0,m]$:
\begin{eqnarray*}
g(m)-g(x)=\int_x^mg'(u)du\geq (m-x)g'(m).
\end{eqnarray*}
Since $g$ is strictly concave and $g'(\infty)=0$, $g'(m)>0$ and we have for $m\geq x$:
\begin{eqnarray*}
\frac{g(m)}{g'(m)}-m\geq\frac{g(x)}{g'(m)}-x.
\end{eqnarray*}
By (\ref{hypginfty}) and monotonicity of $g$, for $x>m_0$, $g(x)>0$ and $\lim_{m\to+\infty}\frac{g(x)}{g'(m)}-x=+\infty$ so that $\lim_{m\to+\infty}\frac{g(m)}{g'(m)}-m=+\infty$.
Since $G'(m)=-\frac{gg''(m)}{(g'(m))^2}$ has the same sign as $g(m)$ by strict concavity of $g$, one easily concludes.
\ep
\\ \\
\noindent{\bf Proof of Theorem \ref{propnash1}}
The characterization conditions \eqref{nash11} and \eqref{nash12} for a Nash equilibrium follow from
the optimal expressions (\ref{opt2}) and Proposition \ref{PropOptimumFirm1bis}. 
Thus, the only thing to check is the existence of an infinite number of solutions in ${\mathbb R}_+\times{\mathbb R}\times{\mathbb R}_+\times {\mathbb R}_+\times{\mathbb R}_+\times {\mathbb R}_+$ to these equations.
\\
By \eqref{minoc} and since $(1+r_J)(D+\hat a)+\hat e >0$, 
\begin{equation*}
 C(\hat \beta)+(1+r_J)(D+\hat a)+\hat e-\frac{u}{u+v}\delta\psi(\hat a)>0.
  \end{equation*}
  
Therefore there exists infinitely many couples  $(\alpha,m)\in \mathbb{R}_+^2$ such that
$$\alpha+G(m)=C(\hat \beta)+(1+r_J)(D+\hat a)+\hat e-\frac{u}{u+v}\delta\psi(\hat a)$$
namely the couples $\left(C(\hat \beta)+(1+r_J)(D+\hat a)+\hat e -\frac{u}{u+v}\delta\psi(\hat a)-G(x),x\right)$ where $x\in[0,G^{-1}(C(\hat \beta)+  (1+r_J)(D+\hat a)+\hat e-\frac{u}{u+v}\delta\psi(\hat a))]$.
\ep

\subsection{Stackelberg equilibrium in Situation 1, firm $I$ is  leader}
\label{Stac1}
As a preliminary, we prove Proposition \ref{eq5.5}, useful for Stackelberg equilibria in both
situations.
\begin{proposition}
 \label{eq5.5}
Equation (\ref{eqbetaStackelberg}) admits at least one solution 
 $\hat\beta.$ Moreover, all solutions are smaller than $\frac{v}{u+v}$.
 \end{proposition}
\proof
Let us recall (\ref{eqbetaStackelberg}): $h(u\beta)-h(v(1-\beta))={S} (\beta)$.
By Lemma \ref{lem:h}, as $\beta$ goes from $-\infty$ to $+\infty$, the left-hand side of this equation is increasing from $h(-\infty)-h(+\infty)<0$ to 
 $h(+\infty)-h(-\infty)>0$ and is null for $\beta=\frac{v}{u+v}.$
 \\
 For $\beta<1$, since $\varphi'(e^*(\beta))=\frac{1}{1-\beta}$, we have
 $S(\beta)=(\varphi'( e^*(\beta))+(b^e)'( e^*(\beta))-1)\frac{(\varphi'( e^*(\beta)))^2}{\varphi''( e^*(\beta))}$
 so that,  by concavity of $\varphi,$ the sign of $S(\beta)$ is equal to the one of
\\
  $1-\varphi'( e^*(\beta))-(b^e)'( e^*(\beta))$. 
 Remember that when $\beta$ goes from $-\infty$ to $1$, $e^*(\beta)$ is decreasing from $+\infty$ to $0$,
  $\varphi$ and $b^e$ are concave, so
 $\beta\mapsto 1-\varphi'( e^*(\beta))-(b^e)'( e^*(\beta))$ is decreasing, from $1$ ($\varphi'(+\infty)=(b^e)'(+\infty)=0$) to $-\infty$ ($\varphi'(0)=(b^e)'(0)=+\infty$). Since $\varphi'+(b^e)'$ is monotonic, there exists a unique $\beta_0$ such that $\varphi'(\beta_0)+(b^e)'(\beta_0)=1$, so
 $\beta<\beta_0\Rightarrow S(\beta)>0,$ $\beta>\beta_0\Rightarrow S(\beta)<0.$
 As a consequence, there exists a solution $\hat\beta$ to (\ref{eqbetaStackelberg}).
For $\beta\geq 0$, $S(\beta)$ is negative since in (\ref{def:S}), the numerator is positive whereas the denominator is negative by concavity of $\varphi$. Hence $\beta_0<0$. Moreover $S(\frac{v}{u+v})<0,$ so any solution $\hat\beta$ belongs to $(\beta_0,\frac{v}{u+v}).$
\ep
\\

\noindent
{\bf Proof of Theorem \ref{PropStackelbergpb1}}
\\
 Firm $I$ has to find $(\alpha,\beta,\gamma)$ maximising 
$$\hspace*{-0.5cm}  b^a(a^*(\beta))- \mathbb{E} \left[ e^{-u  \big(b^m(m^*(\gamma))+b^e(e^*(\beta))-\alpha-\beta(\mu~-\varphi(e^*(\beta))-\delta\psi(a^*(\beta)))-\gamma{g}(m^*(\gamma))\big) }\right].
$$
Since the inverse function of $m^*(\gamma)=(g')^{-1}(1/\gamma)$ is the increasing bijection $\gamma^*(m)=\frac{1}{g'(m)}$ from $\mathbb{R}^+$ onto itself, the maximisers are the triplets $(\alpha,\beta,\gamma^*(m))$ with $(\alpha,\beta,m)$ maximising
$$\tilde I^1(\alpha,\beta,m):=
b^a(a^*(\beta))- \mathbb{E} \left[ e^{-u  \big(b^m(m)+b^e(e^*(\beta))-\alpha-\beta(\mu~-\varphi(e^*(\beta))-\delta\psi(a^*(\beta)))-\frac{g}{g'}(m)\big)  }\right]
$$
under the constraint
\begin{eqnarray*}
&\hspace*{-10cm}-\tilde J^1(\alpha,\beta,m):=\\
&\mathbb{E}  \left[ e^{-v (\alpha+(\beta-1)(\mu~-\varphi(e^*(\beta))-\delta\psi(a^*(\beta)))-e^*(\beta)+\frac{g}{g'}(m)-m -(1+r_J)(D+a^*(\beta) )) }\right] \leq 1,
 \end{eqnarray*}
where, by a slight abuse of notations, we still denote by $\tilde I^1$ and $\tilde J^1$ the functions obtained by applying the change of variable $(\alpha,\beta,\gamma)\to(\alpha,\beta,m)$ to the ones introduced in Section \ref{stackincominf}.
We also recall  the application   $C:\mathbb{R}\to\mathbb{R}$ defined in (\ref{defA})
$$C(\beta)=\frac{1}{v}\ln\mathbb{E}\left[e^{v(1-\beta)(\mu~-\varphi\circ e^*(\beta))}\right].$$

\noindent Setting   
 \begin{eqnarray}
 &{\cal A}:=\Big\{(\alpha,\beta,m)\in \mathbb{R}_+\times\mathbb{R}\times\mathbb{R}_+;~-\tilde J^1(\alpha,\beta,m)\leq 1 \},
 \end{eqnarray}
the optimization problem for firm $I$ then writes:
\begin{eqnarray*}
\sup_{(\alpha,\beta,m)\in{\cal A}}&\tilde I^1(\alpha,\beta,m).
\end{eqnarray*}
We will prove the existence of a maximizer for this problem, and therefore of a Stackelberg equilibrium, by checking that we can restrict the set ${\cal A}$ to a compact subset. Notice first that ${\cal A}\neq\emptyset$. In fact, one can easily check that for any  $\beta\in\mathbb{R}$ and $m\geq 0$, one can choose $\alpha$ large enough so that $(\alpha,\beta,m)\in{\cal A}$.
\\

\begin{proposition}
\label{Lem4.1}
We have $\sup_{(\alpha,m)\in\mathbb{R}_+^2}\tilde I^1(\alpha,\beta,m)\to-\infty$ when $\beta\to\infty$.
Moreover, there exists $\bar\beta\in\mathbb{R}$, not depending on $v>0$, such that the supremum over ${\cal A}$ is attained if and only if the supremum over ${\cal A}\cap\{\beta \in (-\infty,\bar\beta]\}$ is attained, and both supremum are equal.
 \end{proposition}
\proof
To prove this proposition, we need the following result which applies to functions $\varphi$ and $\psi$ :
\begin{lemma}\label{Lemma_equivalentphiprime}
For any increasing, strictly concave $C^1$ and bounded function $\phi$, \\
$\phi'(x)=\circ(1/x)$ when $x\to\infty$ and $y(\phi')^{-1}(y)\to 0$ when $y\to 0$.
\end{lemma}
\proof
 Integrating by parts, we get for $x\geq 1$, 
 \begin{equation}
 \label{ippphi}
 \int_1^xu\phi''(du)=x\phi'(x)-\phi'(1)-\phi(x)+\phi(1)
\end{equation} 
where $\phi''(du)$ denotes the negative measure equal to the 
second order distribution derivative of $\phi$. Since $\phi$ is increasing and concave, 
the terms $-x\phi'(x)$ and $\int_1^xu\phi''(du)$ are non-positive on $[1,+\infty)$. 
The boundedness of $\phi$ then implies their boundedness on $[1,+\infty)$.
 Since $\int_1^xu\phi''(du)$ and $\phi(x)$ are monotonic and bounded, they admit
  finite limits when $x\to\infty$. By (\ref{ippphi}), $x\phi'(x)$ admits a finite
   limit as well, denoted $\ell$. Since $\phi$ is bounded, $\phi'$ is integrable 
   on $[1,+\infty)$, which implies $\ell=0$ and gives the result. Let $x=(\phi')^{-1}(y)$.
    When $y\to 0$, $x$ goes to $\infty$ and $y(\phi')^{-1}(y)= \phi'(x)x$ goes to $0$. 
    \\
 \ep
 
Since $b^m$ is such that ${b^m}'(\infty)=0$, $\lim_{m\to+\infty}  m-b^m(m)  =+\infty$ thus, using the first assertion in Lemma \ref{Lemmalimginfty},
\begin{equation}\label{eqm}
\lim_{m\to+\infty}  \frac{g}{g'}(m)-b^m(m)  =+\infty.
\end{equation}
 Hence $e^{u\big(\frac{g}{g'}(m)-b^m(m)\big)}$ goes to infinity when $m\to\infty$ and there exists a constant $c>0$ not depending on $v>0$ such that for any $m\geq 0$, $e^{u\big(\frac{g}{g'}(m)-b^m(m)\big)}\geq c$.
  \\
  For any $\beta\geq 1$, $e^*(\beta)=0$ and since $\tilde I^1$ is decreasing with respect to $\alpha$, we have for $(\alpha,m)\in\mathbb{R}_+^2$,
\begin{eqnarray}
\label{betamax}
\tilde I^1(\alpha,\beta,m)&\leq& \tilde I^1(0,\beta,m) \nonumber \\ 
&\leq& b^a(a^*(\beta))  -ce^{-ub^e(0)}\;\mathbb{E}\left[e^{u\beta(\mu~-\varphi(0)-\delta \psi(a^*(\beta)))}\right]\nonumber \\ 
&\to& -\infty~\mbox{when}~\beta\to\infty.
\end{eqnarray}
Indeed, if $\delta>0$, then $a^*(\beta)=(\psi')^{-1}\left(\frac{1+r_J}{(\delta(1-\beta))^+}\right)=0$ for $\beta\geq 1$ and the result is obvious. 
\\
Otherwise for $\delta\leq 0,$ $b^a(a^*(\beta))=o(\beta)$  when $\beta \to\infty$ (indeed $b^a(x)=o(x)$ when $x\to\infty$ and $a^*(\beta)=o(\beta)$ when $\beta \to\infty$, see Lemma \ref{Lemma_equivalentphiprime}),  $\mu~-\varphi(0)-\delta \psi(a^*(\beta)) > \mu~-\varphi(+\infty)-\delta^+ \psi(+\infty)>0 $. Therefore $e^{u\beta(\mu~-\varphi(0)-\delta \psi(a^*(\beta)))}$ goes to $+\infty$ faster than $b^a(a^*(\beta))$ and, by Fatou Lemma,
$$\liminf_{\beta \to\infty}\mathbb{E}\left[e^{u\beta(\mu~-\varphi(0)-\delta \psi(a^*(\beta)))-\ln(b^a(a^*(\beta)))}\right]=\infty$$
so that (\ref{betamax}) holds.
\\
 Since $((1+r_J)D,1,0)\in{\cal A}$ for any $v>0$, one deduces the existence of $\bar \beta\in\mathbb{R}$ not depending on $v$ such that 
$$\sup_{(\alpha,\beta,m)\in \mathbb{R}_+\times (\bar \beta,+\infty)\times \mathbb{R}_+}\tilde I^1(\alpha,\beta,m)<\sup_{(\alpha,\beta,m)\in{\cal A}}\tilde I^1(\alpha,\beta,m)$$
and the supremum over ${\cal A}$ is attained if and only if the supremum over
\\
 ${\cal A}\cap\{\beta \in (-\infty,\bar\beta]\}$ is attained.  Moreover, if the suprema are attained,  they are equal.
\ep
\\

 We now decompose the optimisation on ${\cal A}\cap\{\beta \in (-\infty,\bar\beta]\}$ according to the positivity of $\alpha$.
\begin{lemma}
\label{L4.2}
Let 
\begin{eqnarray}
 &{\cal B}:=\Big\{(\alpha,\beta,m)\in \mathbb{R}_+\times(-\infty,\bar\beta]\times\mathbb{R}_+;~-\tilde J^1(\alpha,\beta,m)=1\Big\},
 \label{calBB}
 \\
& {\cal C}:=\Big\{(0,\beta,m);~\beta\in(-\infty,\bar\beta],~m\geq 0,~-\tilde J^1(0,\beta,m)\leq 1\Big\}.
\label{calCC}
\end{eqnarray}
The supremum over ${\cal A}$ is attained if and only if the supremum over ${\cal B}\cup{\cal C}$ is attained.
\end{lemma}
\proof 
Since $\tilde I^1$ is decreasing with respect to $\alpha$ and $\tilde J^1$ is continuous, if $\alpha>0$ and 
\\
$-\tilde J^1(\alpha,\beta,m)<1$, then there exists $\eps>0$ such that $\tilde I^1(\alpha-\eps,\beta,m)>\tilde I^1(\alpha,\beta,m)$, while
 $-\tilde J^1(\alpha-\eps,\beta,m)\leq 1$. Therefore the supremum over ${\cal A}\cap\{\alpha>0\}$ is the same as the supremum over $$\tilde{\cal B}:=\Big\{(\alpha,\beta,m)\in (0,+\infty)\times(-\infty,\bar\beta]\times\mathbb{R}_+;~\tilde J^1(\alpha,\beta,m)=-1\Big\}.$$
 Noticing then that ${\cal A}$ is closed and that ${\cal B}:=Cl(\tilde{\cal B})$, we easily conclude.
\ep\\

\noindent Lemmas \ref{L4.3} and \ref{L4.4} are devoted to the supremum over $\mathcal{B}$. In Lemma  \ref{L4.3}, we check that the supremum is attained on a compact subset $\mathcal{K}$ not depending on $v>0$. In Lemma \ref{L4.4}, we derive optimality  conditions satisfied by a maximizer with positive $\alpha$. 
\begin{lemma}
\label{L4.3}
The function $k$
defined in (\ref{defk}):
$$k:\beta\to e^{u(Id-b^e)\circ e^*(\beta)}   e^{-u \delta \psi \circ a^*(\beta)}e^{u(1+r_J)(D+a^*(\beta))}   e^{uC(\beta)}\mathbb{E}\left[e^{u\beta\big(\mu-\varphi\circ e^*(\beta)\big)}\right]$$
 depends on $v>0$ through the function $C(\beta)$ defined in (\ref{defA}), goes to $\infty$ uniformly in $v>0$ when $\beta\to-\infty$.
 \\
 Moreover, in the set ${\cal B},$  $\alpha$ is a 
 continuous function of $(\beta,m)$ and there exists a compact set ${\cal K}\subset{\cal B}$
 not depending on $v>0$ such that $\sup_{{\cal B}}\tilde I^1=\sup_{{\cal K}}\tilde I^1$. In particular 
 the supremum on ${\cal B}$ is attained.
\end{lemma}
\proof
 Let us consider $\sup_{{\cal B}}\tilde I^1$. Notice that $C(\beta)$ given by (\ref{defA}) depends on $v>0$. Since $\tilde J^1(\alpha,\beta,m)=-1$ on ${\cal B}$, we know that
 $$\alpha=\alpha(\beta,m):=C(\beta)-(1-\beta)\delta\psi(a^*(\beta))+e^*(\beta)+m-\frac{g}{g'}(m)+(1+r_J)(D+a^*(\beta)),$$ so that
\begin{equation}
   \tilde I^1(\alpha(\beta,m),\beta,m)=b^a(a^*(\beta))-e^{u(Id-b^m)(m)} k(\beta).\label{i1albetm}
\end{equation}
In $\tilde I^1(\alpha(\beta,m),\beta,m)$, we have $b^a(a^*(\beta)) $ minus the product of two positive functions, the first one depending only on $m$, the second one, $k$, only on $\beta$. Since $(b^m)'(\infty)=0$,  $e^{u(Id-b^m)(m)}\to\infty$ when $m\to\infty$, therefore the first function is bounded from below by a positive constant $c$ and goes to infinity when $m$ goes to infinity. 
\\
 Let us now examine the function $k$. By Jensen's inequality, $\E\left[e^{u\beta(\mu-\varphi(e^*(\beta)))}\right]\geq e^{u\beta(\E[\mu]-\varphi(e^*(\beta)))}$; then, using \eqref{minC} and $e^{u(1+r_J)(D+a^*(\beta))}\geq 1$ it yields
$$
k(\beta)\geq e^{u(Id-b^e)(e^*(\beta))}
e^{u\Big(\mathbb{E}[\mu]~-\varphi(e^*(\beta))-\delta \psi(a^*(\beta))  \Big)}.
$$
\noindent 
Using $\mathbb{E}\left[\mu~\right]\geq \varphi(x)+\delta^+ \psi(y)\mbox{ for any }x,y$, we deduce that
 $k(\beta)\geq e^{u(Id-b^e)\circ(\varphi')^{-1}(\frac{1}{(1-\beta)^+})}$ where the right-hand side does not depend on $v>0$ and goes to infinity when $\beta\mapsto-\infty$
since  $b^e(x)=o(x)$   when $x\to\infty$ (recall that $(b^e)'(\infty)=0$). Finally
using   $b^a(x)=o(x)$ and $(\varphi')^{-1}(x) \sim (\psi')^{-1}(x) $ for $x \rightarrow 0$, we conclude when $\beta\to -\infty$:  
$$ \tilde I^1(\alpha,\beta,m)
\leq b^a\circ(\psi')^{-1}(\frac{1+r_J}{\delta^+(1-\beta)})   -c   e^{u(Id-b^e\\)\circ(\varphi')^{-1}(\frac{1}{1-\beta})}  \\
\to-\infty.$$
Since $\alpha(1,0)=(1+r_J)D>0$ and $\tilde I^1(\alpha(1,0),1,0)$ do not depend on $v>0$ and $(\alpha(1,0),1,0)\in{\cal B}$ for any $v>0$, 
we deduce the existence of $\underline{\beta}>-\infty$ and $\bar{m}$ in $[0,+\infty)$ not depending on $v$ such that $\sup_{\cal B}\tilde I^1=\sup_{{\cal B}\cap\{\beta\in[\underline{\beta},\bar\beta],m\in[0,\bar{m}]\}}\tilde I^1$. 

To conclude that $\sup_{\cal B}\tilde I^1$ is attained on a compact ${\cal K}\subset {\cal B}$ not depending on $v$, one still
has to make sure that $\alpha$ remains bounded when $v$ varies. Since for any $\alpha$,
\begin{align*}
 &\sup_{\beta\in[\underline{\beta},\bar\beta],m\in[0,\bar{m}]}
 \tilde I^1(\alpha,\beta,m)\leq\sup_{\beta\in[\underline{\beta},\bar\beta]}
 b^a(a^*(\beta))
   \\ &-e^{u\alpha}\inf_{\beta\in[\underline{\beta},\bar\beta],m\in[0,\bar{m}]}
   \mathbb{E} \left[ e^{-u  \big(b^m(m)+b^e(e^*(\beta))-\beta(\mu~-\varphi(e^*(\beta))-\delta\psi(a^*(\beta)))-\frac{g}{g'}(m)\big)  }\right]
\end{align*}
where the right-hand side tends to $-\infty$ as $\alpha\to\infty$, one concludes that there exists a compact set ${\cal K}\subset{\cal B}$ not depending on $v>0$ such that $\sup_{\cal B}\tilde I^1=\sup_{\cal K}\tilde I^1$.
\ep
\\
\begin{lemma}
\label{L4.4}
Let $(\hat\alpha,\hat\beta,\hat m)$ with $\hat\alpha>0$ be such that the maximum on ${\cal B}$ is attained at $(\hat\alpha,\hat\beta,\hat m)$. Then necessarily 
\begin{equation}
\label{eqdegamma'}
\hat m=[(b^m)']^{-1}(1)>0,
\end{equation}
and there are two different cases, depending on the sign of $\delta$ :\\
$\bullet$ If $\delta>0$,  
$$ h(u\hat{\beta})-h(v(1-\hat{\beta}))=S(\hat{\beta})+\widetilde{S} (\hat{\beta}), ~\mbox{and }~ \hat{\beta}<1. \, \, \,    (\ref{eqbetaStackelbergpb1positif}) $$
$\bullet$ If $\delta<0$,
$$ h(u\hat{\beta})-h(v(1-\hat{\beta}))=S(\hat{\beta}) ,  ~\mbox{and}~ \hat{\beta}<1, \, \, \,  (\ref{eqbetaStackelberg})$$
or
$$
h(u\hat{\beta})-h(v(1-\hat{\beta}))=\widetilde{S} (\hat{\beta}),  ~\mbox{and}~ \hat{\beta} \geq 1. \, \, \,   (\ref{eqbetaStackelbergpb1negatif})
$$
with $h$, $S$ and $\widetilde{S}$ respectively defined in (\ref{def:h}), (\ref{def:S}) and (\ref{def:Stilde}).
\end{lemma}
\proof
Let $(\hat\alpha,\hat\beta,\hat m)$ be such that the maximum on ${\cal B}$ is attained at
 $(\hat\alpha,\hat\beta,\hat m)$. Since $(b^m)'(0)=+\infty$, $m\mapsto m-b^m(m)$ is decreasing
  in a neighborhood of $0$, so that from (\ref{i1albetm}), $\hat m>0$. Assume moreover
   that $\hat\alpha>0$. Then the mapping $(\beta,m)\mapsto \tilde I^1(\alpha(\beta,m),\beta,m)$ where $\alpha(\beta,m)$ is defined at the beginning of the proof of Lemma \ref{L4.3}
   admits a local maximum at $(\hat\beta,\hat m)$ and therefore the first order conditions
    are satisfied (notice that $\bar{\beta}$ may be increased), i.e.
     $\frac{\partial}{\partial \beta} \tilde I^1(\alpha(\hat\beta,\hat m),\hat\beta,\hat m)=\frac{\partial}{\partial m} \tilde I^1(\alpha(\hat\beta,\hat m),\hat\beta,\hat m)=0$. 
     The second one gives
\begin{equation}\label{eqdegammapb1}
\hat m=[(b^m)']^{-1}(1)>0.
\end{equation}
The computation of $\hat\beta$ is more tricky and depends on the coefficient $\delta$.

\vspace{5mm}

1) \underline{ $\delta>0$}: 
$\frac{\partial}{\partial\beta}
\tilde I^1(\alpha(\hat\beta,\hat m),\hat\beta,\hat m)   =
$
\begin{align*}
 \begin{cases}
\bullet -e^{u(Id-b^m)(\hat m)}
 \mathbb{E}\left[u\nu
e^{u[\beta\mu+C(\beta)]}\right]~\mbox{ if }~\beta\geq 1
\\
\bullet (b^a)'\circ a^*(\beta)
\frac{1+r_J}{\delta(1-\beta)^2 \psi''\circ a^*(\beta)}
\\
-u\mathbb{E}\Big(\frac{\beta(1+r_J)}{\beta-1}(a^*)'(\beta)-
\frac{\beta}{1-\beta}(e^*)'(\beta)-(b^e)'\circ e^*(\beta)(e^*)'(\beta)+\nu\Big)
e^{u(Id-b^m)(\hat m)}
\\
e^{u[(Id-b^e)\circ e^*(\beta)- \delta \psi\circ a^*(\beta)
 +\beta(\mu-\varphi\circ e^*(\beta))+(1+r_J)(D+a^*(\beta))+C(\beta)]}~
 \mbox{if}~\beta<1
\end{cases}
\end{align*}
where $\nu:=\mu-h(v(1-\beta))$. Then we get
\begin{align*}
\frac{\partial}{\partial \beta} \tilde I^1(\alpha(\hat\beta,\hat m),\hat\beta,\hat m)   = 0  ~\mbox{iff}~ 
\begin{cases}
h(u\beta)-h(v(1-\beta))=0   ~\mbox{if}~\beta\geq 1,
\\
h(u\beta)-h(v(1-\beta))=S(\beta) +\widetilde{S} (\beta) ~\mbox{ if}~\beta<1,
\end{cases}
\end{align*}
recalling
$$S(\beta):= \frac{\frac{\beta}{1-\beta}+(b^e)'\circ(\varphi')^{-1}\Big(\frac{1}{1-\beta}\Big)}{(1-\beta)^2\varphi''\circ(\varphi')^{-1}\Big(\frac{1}{1-\beta}\Big)}.
$$

$$
\widetilde{S}(\beta):= \frac{1+r_J}{\delta(1-\beta)^2 (\psi'')(a(\beta)) } \left(  (1+r_J)\frac{\beta}{1-\beta} +  \frac{(b^a)'(a(\beta))}{ u   e^{u(Id-b^m)(\hat m)}  k(\beta)  }      \right)
$$
As seen in Lemma \ref{lem:h}, for $\beta\geq 1$, $h(u\beta)-h(v(1-\beta))>0$, thus 
$\hat\beta<1$ and we study the equation 
(\ref{eqbetaStackelbergpb1positif})
$$
h(u\beta)-h(v(1-\beta))=S(\beta)+\widetilde{S} (\beta),  ~\mbox{for}~ \beta<1. 
$$
The left hand side is positive  for $\beta>\frac{v}{u+v} $. The functions $S$ and $\widetilde{S}$ are negative on $[0,1]$ thus $\hat\beta<\frac{v}{u+v}$.

\vspace{5mm}

2) \underline{ $\delta<0$}: $ \frac{\partial}{\partial \beta }\tilde I^1(\alpha(\hat\beta,\hat m),\hat\beta,\hat m)   =$
\begin{align*}
\begin{cases}
\bullet(b^a)'\circ a^*(\beta)
\frac{1+r_J}{\delta(1-\beta)^2 \psi''\circ a^*(\beta)}
-u\mathbb{E}\left[(\frac{\beta(1+r_J)}{\beta-1}
(a^*)'(\beta)+\nu)\right.
\\
\left.\times e^{u[(Id-b^m)(\hat m)-b^e(0)- \delta \psi\circ a^*(\beta) 
+\beta(\mu-\varphi(0))+(1+r_J)(D+a^*(\beta))+C(\beta)]}\right]~\mbox{if}~\beta\geq 1
\\
\\
\bullet-u\mathbb{E}\left[\big(-\frac{\beta}{1-\beta}(e^*)'(\beta)-
(b^e)'\circ e^*(\beta)(e^*)'(\beta)+\nu\big)\right.
\\
\left.\times e^{u[(Id-b^m)(\hat m)+(Id-b^e)
\circ e^*(\beta)-\delta\psi(0)+\beta(\mu-\varphi
\circ e^*(\beta))+(1+r_J)D+C(\beta)]}\right]~\mbox{if}~\beta<1.
\end{cases}
\end{align*}

Thus, 
\begin{align*}
\frac{\partial}{\partial \beta}\tilde I^1(\alpha(\hat\beta,\hat m),\hat\beta,\hat m)   = 0 ~\mbox{iff}~ \begin{cases}
h(u\beta)-h(v(1-\beta))=  \widetilde{S} (\beta) ~\mbox{if}~\beta\geq 1
\\
h(u\beta)-h(v(1-\beta))= S (\beta)  ~\mbox{if}~\beta<1.
\end{cases}
\end{align*}

Thus the optimal $\beta$ on ${\cal B}$ is either the solution of 
equation (\ref{eqbetaStackelberg}) (which is less than $\frac{v}{u+v}$) or the solution of 
$$h(u\beta)-h(v(1-\beta))=\widetilde{S} (\beta),  ~\mbox{for}~ \beta \geq 1. 
$$
\ep


\vspace{0.6mm}

\noindent The next lemma is devoted to the optimisation over $\mathcal{C}$.
\begin{lemma}
\label{L4.5}
The supremum of $(\alpha,\beta,m)\mapsto \tilde I^1(\alpha,\beta,m)$ on ${\cal C}$ defined in (\ref{calCC}) is attained.
\end{lemma}
\proof
Let us then consider $\sup_{{\cal C}}\tilde I^1$. We have $\alpha=0$ on ${\cal C}$. Since $(b^m)'(\infty)=0$ and, by Lemma \ref{Lemmalimginfty}, $G(m)=\frac{g}{g'}(m)-m$ is increasing on $[m_0,+\infty)$, the mapping $m\mapsto b^m-\frac{g(m)}{g'(m)}$ is decreasing for sufficiently large $m$. Thus there is a constant $\bar m\in [m_0,+\infty)$ not depending on $v>0$ such that for any $\beta\in\mathbb{R}$, $m\mapsto \tilde I^1(0,\beta,m)$ is decreasing for $m\geq \bar m$. Therefore, writing
$${\cal C}_1:={\cal C}\cap\{m\in[0,\bar m]\} \quad\mbox{and}\quad {\cal C}_2:={\cal C}\cap\{m\geq\bar m~\mbox{and}~\tilde J^1(0,\beta,m)=-1\},$$
the supremum of $\tilde I^1$ over ${\cal C}$ is attained iff the supremum of $\tilde I^1$ over ${\cal C}_1\cup{\cal C}_2$ is attained.\\

(i) We prove that ${\cal C}_1$ is compact. 
The condition $-\tilde J^1(0,\beta,m)\leq 1$ is equivalent to
$$C(\beta)-(1-\beta)\delta\psi(a^*(\beta))+e^*(\beta)+m-\frac{g}{g'}(m)+(1+r_J)(D+a^*(\beta))\leq 0.$$
If $\beta<1$, by \eqref{minoc}, it implies $e^*(\beta)+m-\frac{g}{g'}(m)+(1+r_J)(D+a^*(\beta))\leq 0$. Since $\lim_{\beta\to-\infty}e^*(\beta)=+\infty$
 and $a^*(\beta)\geq 0$, while $m\in[0,\bar{m}]$, this implies that there exists $\underline\beta\in\mathbb{R}$ not depending on $v>0$, such that ${\cal C}_1\subset \{0\}\times[\underline\beta,\bar\beta]\times[0,\bar m]$, $\bar \beta$ being defined in Proposition~\ref{Lem4.1}. Since ${\cal C}_1$ is closed by continuity of $\tilde J^1$, it is compact and therefore the supremum over ${\cal C}_1$ is attained.\\

(ii) On ${\cal C}_2$, we have $$C(\beta)-(1-\beta)\delta\psi(a^*(\beta))+e^*(\beta)+(1+r_J)(D+a^*(\beta))=\frac{g(m)}{g'(m)}-m=G(m).$$ 
Since $m\geq \bar{m}\geq m_0$, using Lemma \ref{Lemmalimginfty} and the inverse function $G^{-1}$, $m=H(\beta):=G^{-1}\big[C(\beta)-(1-\beta)\delta\psi(a^*(\beta))+e^*(\beta)+(1+r_J)(D+a^*(\beta))\big]$.
Thus we have $\sup_{{\cal C}_2}\tilde I^1=\sup_{{\cal C}_2\cap\{\beta \in (-\infty,\bar \beta]\}}\tilde I^1\big(0,\beta,H(\beta)\big)$. \\
Let us now prove that $\tilde I^1(0,\beta,H(\beta))\to-\infty$ uniformly in $v>0$ when $\beta\to-\infty$.
When $\beta\to-\infty$, $e^*(\beta)\to\infty$, while $a^*(\beta)\geq 0$, $C(\beta)-(1-\beta)\delta\psi(a^*(\beta))\geq 0$ by \eqref{minoc}, so that
$C(\beta)-(1-\beta)\delta\psi(a^*(\beta))+e^*(\beta)+(1+r_J)(D+a^*(\beta))\to\infty$
 uniformly in $v>0$, and therefore $m=H(\beta)\to\infty$.
We recall that  $b^m(m)=\circ(m)$, $a^*(\beta)=\circ(-\beta)$ and $e^*(\beta)=\circ(-\beta)$ (see Lemma \ref{Lemma_equivalentphiprime}) and 
$b ^e\circ e^*(\beta)=\circ(-\beta)$,~ when
$\beta\to-\infty$. Using $$\frac{g\circ H(\beta)}{g'\circ H(\beta)}-H(\beta)=C(\beta)-(1-\beta)\delta\psi(a^*(\beta))+e^*(\beta)+
(1+r_J)(D+a^*(\beta)),$$
we deduce that for all $v>0$,
\begin{align}
   -&u\left(b^m\circ H(\beta)+b^e(e^*(\beta))-\beta(\mu~-\varphi(e^*(\beta))-\delta \psi(a^*(\beta)))-\frac{g\circ H(\beta)}{g'\circ H(\beta)}\right)\notag\\
=&-u\big(b^m\circ H(\beta)-H(\beta)+(b^e-Id)(e^*(\beta))-\beta(\mu~-\varphi(e^*(\beta))-\delta \psi(a^*(\beta)))\\&\phantom{-u\big(}-C(\beta)+(1-\beta)\delta\psi(a^*(\beta))-(1+r_J)(D+a^*(\beta))\big)\notag
\\
&\geq-u\beta\left(\mathbb{E}[\mu]~ -\mu~\right)+o(-\beta),
 \mbox{ as }-\beta\to\infty,\label{inegaliteenbetapb1}
\end{align}
where the last inequality follows from (\ref{minC}) ,
$o(-\beta)$ being uniform in $v>0$.
\\
Since $\mu$ is not $\pr$ constant, 
\begin{equation}
\label{controlmu}
\exists \eps>0, \exists{\cal D}\subset \Omega ~ \mbox{  s.t. }~ \pr({\cal D})\geq \eps\mbox{ and  }\mu+\eps\leq  \mathbb{E}[\mu]~  \mbox{ on }{\cal D}.
\end{equation}
 Using (\ref{inegaliteenbetapb1}), we compute for $-\beta$ sufficiently large not depending on $v>0$ :
$$
\tilde I^1(0,\beta,H(\beta))\leq  b^a(a^*(\beta))-\mathbb{E}[ 1_{{\cal D}}e^{-\frac{u\beta}{2}\left(\mathbb{E}[\mu]~-\mu\right)}]
\leq b^a(a^*(\beta)) -\eps e^{-\frac{u\beta\eps}{2}}.$$
Since $b^a(a^*(\beta))=o(-\beta)$ when $\beta\to-\infty $, the right-hand side goes to $-\infty$ uniformly in $v>0$ when $\beta\to-\infty$. 
Since $(0,1,G^{-1}((1+r_J)D))\in{\cal C}$ for any $v>0$ and $\tilde I^1(0,1,G^{-1}((1+r_J)D))$ does not depend on $v>0$, one deduces the existence of $\underline\beta>-\infty$ such that if $\sup_{\cal C}\tilde I^1=\sup_{{\cal C}_2}\tilde I^1$ then $\sup_{\cal C}\tilde I^1=\sup_{{\cal C}_2\cap\{\beta\in[\underline\beta,\bar\beta]\}}\tilde I^1$. Now 
 \begin{align*}
   &\sup_{\beta\in[\underline{\beta},\bar\beta]}\tilde I^1(0,\beta,m)\leq\sup_{\beta\in[\underline{\beta},\bar\beta]}b^a(a^*(\beta))\\&-e^{u\left(\frac{g}{g'}(m)-b^m(m)\right)}\inf_{\beta\in[\underline{\beta},\bar\beta]}\mathbb{E}\left[e^{-u\big(b^e(e^*(\beta))-\beta(\mu~-\varphi(e^*(\beta))-\delta\psi(a^*(\beta)))\big)}\right]
\end{align*}
where, by \eqref{eqm}, the right-hand-side tends to $-\infty$ as $m\to\infty$. Hence there exists $\bar{m}\in[0,+\infty)$ not depending on $v>0$ such that if $\sup_{\cal C}\tilde I^1=\sup_{{\cal C}_2}\tilde I^1$ then $\sup_{\cal C}\tilde I^1=\sup_{{\cal C}_2\cap\{\beta\in[\underline\beta,\bar\beta],m\in[0,\bar{m}]\}}\tilde I^1$.
\ep
\\

To conclude the proof of Theorem \ref{PropStackelbergpb1},  Lemmas \ref{L4.2}, \ref{L4.3} and \ref{L4.5} prove that  the maximum of $\tilde I^1$ over ${\cal A}$ is attained at $(\hat \alpha,\hat \beta,\hat m)$ which belongs either to ${\cal B}$ or to ${\cal C}$, and a Stackelberg equilibrium exists. Moreover if $\hat \alpha>0$, then $(\hat \alpha,\hat \beta,\hat m)\in{\cal B}$ and the maximum of $\tilde I^1$ over ${\cal B}$ is attained at $(\hat \alpha,\hat \beta,\hat m)$. The equilibrium 
 characterization given in the statement of Theorem \ref{PropStackelbergpb1} then follows from Lemma \ref{L4.4}. 
 \ep

\section{Proofs in Situation 2}
\label{sec:proof2}

\subsection{Best responses in Situation 2}
\label{ssec:4.2}
%


Let $e\geq 0$ and $m\geq 0$ be given and constant. We recall that $F(a)=b^a(a)-f((1+r_I)(D+a))$
is assumed to be  strictly concave, $F'(0)>0$, possibly infinite, and $F'(\infty)=-\infty$. We introduce 
$$K_a(\alpha,\beta,\gamma):=-\mathbb{E}\left[e^{-u\big(b^m(m)+b^e(e)-\alpha-\beta C^{op}-\gamma g(m)\big)}\right]$$
 Then we get the following optimization problem:
\begin{eqnarray}
&\sup_{a,\alpha,\beta,\gamma}F(a)&+K_a(\alpha,\beta,\gamma)\mbox{ with}~a\geq 0,~\alpha\geq 0,~\gamma\geq 0,\\
&\mbox{and}&~\mathbb{E}\left[e^{-v\big(\alpha+(\beta-1)(\mu~-\varphi(e)-\delta\psi(a))-e+\gamma g(m)-m\big)}\right]\leq 1.\nonumber
\end{eqnarray}

We have the following result:
\begin{proposition}
\label{PropOptimumFirm1}
Let $e\geq 0$ and $m\geq 0$ be given and constant and let $\beta^*:=\frac{v}{u+v}$. 
There exists optimal controls and $(a,\alpha,\beta,\gamma)$ is optimal if and only if it satisfies: 
$\beta=\beta^*$, 
$a\in\arg\max_{a\geq 0} \left[ F(a)-e^{-u\delta\psi(a)}e^{-u(b^m(m)+b^e(e)-e-m)}e^{(u+v)C_e(\beta^*)}\right]$ 
and 
\begin{equation}
\alpha+\gamma g(m)=C_e(\beta^*)-\frac{u}{u+v}\delta\psi(a)+e+m
\label{infinitsol}
\end{equation}
 with $\alpha\geq 0$ and $\gamma\geq 0$ and  $C_e(\beta):=\frac{1}{v}\ln\mathbb{E}\left[e^{v(1-\beta)(\mu-\varphi( e))}\right]$.
\\
Last, if $\delta\geq 0$, then $a^*>0$ and is unique.

\end{proposition}

\begin{remark}
 Notice that $\beta^*\in(0,1)$ and any optimal control $(a,\alpha,\beta,\gamma)$ is such that $\alpha+\gamma>0$.
\end{remark} 

\proof
Let $a\geq 0$ be given for the moment and
$$E_a:=\Big\{(\alpha,\beta,\gamma)\in\mathbb{R}^3;~\mathbb{E}\left[e^{-v\big(\alpha+(\beta-1)C^{op}-e+\gamma g(m)-m\big)}\right]\leq 1\Big\}.$$
In comparison with the definition of the set $E$ introduced in the proof of Proposition \ref{PropOptimumFirm1bis}, the term $(1+r_J)(D+a)$ does not appear in the constraint defining the set $E_a$.
\\
We will first solve, for fixed $a$, the problem of maximization of $K$ on $E_a$, following the same steps as in  Proposition \ref{PropOptimumFirm1bis} :  the optimal $(\alpha^*,\beta^*,\gamma^*)$ for the problem
$$\sup_{(\alpha,\beta,\gamma)\in E_a\cap(\mathbb{R}_+\times\mathbb{R}\times\mathbb{R}_+)}K(\alpha,\beta,\gamma),$$
are exactly the elements of $$\big\{(\alpha,\beta,\gamma)\in\mathbb{R}_+\times\mathbb{R}\times\mathbb{R}_+;~\beta=\frac{v}{u+v},~\alpha+\gamma g(m)=C_e(\beta)-\frac{u}{u+v}\delta\psi(a)+e+m\big\}$$
and the constraint is always binding.
Using  $\beta^*=\frac{v}{u+v}$, we get $\E[e^{u\beta^*(\mu-\varphi(e))}]=e^{vC_e(\beta^*)}$ so that, with the equality $\alpha^*+\gamma^* g(m)=C_e(\beta^*)-\frac{u}{u+v}\delta\psi(a)+e+m$,$$K_a(\alpha^*,\beta^*,\gamma^*)=-e^{-u\delta\psi(a)}e^{-u(b^m(m)+b^e(e)-e-m)}e^{(u+v)C_e(\beta^*)}.$$

Let us then consider  $\sup_{a\geq 0}F(a)+K_a(\alpha^*,\beta^*,\gamma^*)$. Since  $F(a)+K_a(\alpha^*,\beta^*,\gamma^*)\leq F(a)$ and $\lim_{a\to\infty}F(a)=-\infty$, $a\mapsto F(a)+K_a(\alpha^*,\beta^*,\gamma^*)$ attains its maximum on $\mathbb{R}_+$, either at $a^*=0$ or at a point where the first order condition is satisfied.
If $\delta\geq 0$, this function is strictly concave, increasing for small $a$, so that there exists a unique maximum $a^*>0$ that is the unique solution of
$$F'(a)-u\delta \psi'(a)K_a(\alpha^*,\beta^*,\gamma^*)=0.$$\ep

\begin{remark}
If $\delta<0$, then the function $a\mapsto F(a)+K_a(\alpha^*,\beta^*,\gamma^*)$ is decreasing for small $a$, since $\psi'(0)=+\infty$, so that $a=0$ is a local maximum. So the maximum is attained either at  $a^*=0$ or at a solution of $F'(a)+u\delta B_e(\beta^*,  m)\psi'(a)e^{-u\delta\psi(a)}=0$.
\end{remark}

\subsection{Nash equilibrium}

From the previous results, the proof of Theorem \ref{propnash2} follows easily.

\noindent {\bf Proof of Theorem \ref{propnash2}}
\\
The characterization
  conditions (\ref{NashEqui}), (\ref{NashEqui2}) and (\ref{NashEqui3})   follow from Proposition \ref{PropOptimumFirm1} and the optimal expressions (\ref{e*m*}). 
Thus, the only thing to check is that there exists an infinite number of solutions to these equations.  Since $\hat\beta<1$, by \eqref{minoc}, $C(\hat \beta) +\hat e-\frac{u}{u+v}\delta\psi({\hat a})\geq \hat{e}>0$. Therefore there exists infinitely many couples  $(\alpha,m)\in \mathbb{R}_+^2$ such that
$$\alpha+G(m)=C(\hat \beta)++\hat e-\frac{u}{u+v}\delta\psi(\hat a)$$
namely the couples $\left(C(\hat \beta)+\hat e -\frac{u}{u+v}\delta\psi(\hat a)-G(x),x\right)$ where $x\in[0,G^{-1}(C(\hat \beta)+ \hat e-\frac{u}{u+v}\delta\psi(\hat a))]$.
\ep

\subsection{Stackelberg equilibrium in Situation 2, firm $I$ is leader}
\label{stackpb2}
The best response for firm $J$ is given by (\ref{e*m*}), but now the optimization problem for firm $I$ has changed. We recall the continuous mappings $m^*:\mathbb{R}_+\to\mathbb{R}_+$,  $e^*:\mathbb{R}\to\mathbb{R}_+$, $C:\mathbb{R}\to\mathbb{R}$ and $\tilde B:\mathbb{R}\times\mathbb{R}_+\to\mathbb{R}$ defined in (\ref{e*m*}) (\ref{defA}) and (\ref{defB}), 
\begin{eqnarray*}
&e^*(\beta)=(\varphi')^{-1}\left(\frac{1}{(1-\beta)^+}\right),~m^*(\gamma)=(g')^{-1}(1/\gamma).\\
&C(\beta)=\frac{1}{v}\ln\mathbb{E}\left[e^{v(1-\beta)(\mu~-\varphi\circ e^*(\beta))}\right].\end{eqnarray*}

We  are now ready to prove the existence of a Stackelberg equilibrium, as stated in Theorem \ref{PropStackelberg}.\\

\noindent {\bf Proof of Theorem \ref{PropStackelberg}.}\\
Given $\beta$ and $\gamma$, the optimal controls for firm $J$ are given by $e^*$ and $m^*$. 
 Once again, since   $\gamma=\frac{1}{g'(m)}$  yields a bijection between
$m$ and $\gamma$ on $\mathbb{R}^+$,  we only deal with $m.$
Writing:
 \begin{eqnarray*}
 &&\tilde I^2(a,\alpha,\beta,m):=
 \\
 &&F(a)- \mathbb{E}\left[e^{-u\big(b^m(m)+b^e(e^*(\beta))-\alpha-\beta(\mu~-\varphi(e^*(\beta))-\delta \psi(a))-
 \frac{ g}{g'}(m)\big)}\right],
 \\
 &&\tilde J^2(a,\alpha,\beta,m):=
 -\mathbb{E}\left[e^{-v\big(\alpha+
 (\beta-1)(\mu~-\varphi(e^*(\beta))-\delta\psi(a))-e^*(\beta)+
 \frac{g}{g'}(m)-m\big)}\right],
 \\
& &{\cal A}:=\Big\{(a,\alpha,\beta,m)\in\mathbb{R}_+\times\mathbb{R}_+\times\mathbb{R}\times\mathbb{R}_+;~-\tilde J^2(a,\alpha,\beta,m)\leq 1\Big\},
 \end{eqnarray*}
the optimization problem for firm $I$ then writes:
\begin{eqnarray*}
\sup_{(a,\alpha,\beta,m)\in{\cal A}}&\tilde I^2(a,\alpha,\beta,m).
\end{eqnarray*}
We will prove the existence of a maximizer for this problem, and therefore of a Stackelberg equilibrium, by proving that we can restrict the set ${\cal A}$ to a compact subset. Notice first that ${\cal A}\neq\emptyset$. Indeed,  since
 $\tilde J^2(a,\alpha,1,0)=-e^{-v\alpha}$,  for any $a\geq 0$ and $\alpha>0$, $(a,\alpha,1,0)\in{\cal A}$. In fact, one can easily check that for any $a\geq 0$, $\beta\in\mathbb{R}$ and $m\geq 0$, one can choose $\alpha$ large enough so that $(a,\alpha,\beta,m)\in{\cal A}$.
\\

 \noindent The proof will use the following lemmas, very similar to Proposition \ref{Lem4.1} and Lemmas \ref{L4.2},
\ref{L4.3}, \ref{L4.4} and \ref{L4.5}. Nevertheless, we cannot deduce them from
the previous ones because the involved functions are not defined on the same spaces.
\begin{lemma}
\label{L5.1}
There exists $\bar a\in (0,+\infty)$ 
and $\bar\beta\in\mathbb{R}$ not depending on $v$ such that the supremum over ${\cal A}$ is attained if and only if the supremum over $\{(a,\alpha,\beta,m)\in{\cal A}\mbox{ such that }a\in[0,\bar a], \beta\in(-\infty,\bar\beta]\}$ is attained.
\end{lemma}
\proof 
For any $v>0$ and any $(a,\alpha)\in{\mathbb R}_+\times{\mathbb R}_+$, $(a,\alpha,1,0)\in{\cal A}$. We have $\tilde I^2(a,\alpha,\beta,m)\leq F(a)$. Since $\lim_{x\to\infty}F(x)=-\infty$
(cf. the beginning of Subsection \ref{ssec:4.2}) there exists $\bar a>0$ not depending on $v$, such that the supremum over ${\cal A}$ is the same as the supremum over ${\cal A}\cap\{a\in[0,\bar a]\}$. 
Notice that $\tilde I^2$ is decreasing with respect to $\alpha$, so that $\tilde I^2(a,\alpha,\beta,m)\leq \tilde I^2(a,0,\beta,m)$.
Using (\ref{eqm}),  we get
  $e^{u\big(\frac{ g}{g'}(m)-b^m(m)\big)}\to +\infty$ when $m\to +\infty$. Therefore, there exists a constant $c>0$ such that for any $m\geq 0$, $e^{u\big(\frac{ g}{g'}(m)-b^m(m)\big)}\geq c.$
  \\
On the other hand, for any $\beta\geq 1$, since $e^*(\beta)=0$, we compute:
$$
\tilde I^2(a,\alpha,\beta,m)\leq \tilde I^2(a,0,\beta,m)
\leq F(a)-c\;e^{-ub^e(0)}\mathbb{E}\left[e^{u\beta(\mu~-\varphi(0)-\delta \psi(a))}\right]
$$
which goes to $-\infty$ when $\beta\to\infty,$ uniformly in $a\in[0,\bar a].$
\ep

\begin{lemma}
\label{L5.2}
Let
 \begin{eqnarray}
 &{\cal B}:=\Big\{(a,\alpha,\beta,m)\in[0,\bar a]\times\mathbb{R}_+\times(-\infty,\bar\beta]\times\mathbb{R}_+;~\tilde J^2(a,\alpha,\beta,m)=-1\Big\},\label{calB}
 \\
& {\cal C}:=\Big\{(a,0,\beta,m);~a\in[0,\bar a],~\beta\in(-\infty,\bar\beta],~m\geq 0,~\tilde J^2(a,0,\beta,m)\geq -1\Big\}.
\label{calC}
\end{eqnarray}
The supremum over ${\cal A}$ is attained if and only if the supremum over ${\cal B}\cup{\cal C}$ is attained.
\end{lemma}
\proof
Since $\tilde I^2$ is decreasing with respect to $\alpha$ and $\tilde J^2$ is continuous, if $\alpha>0$ and $-\tilde J^2(a,\alpha,\beta,m)<1$, then there exists $\eps>0$ such that $\tilde I^2(a,\alpha-\eps,\beta,m)>\tilde I^2(a,\alpha,\beta,m)$, while $-\tilde J^2(a,\alpha-\eps,\beta,m)\leq 1$. Thus in case of optimum satisfying 
$\alpha>0,$ the constraint is binding.
\ep

\vspace{5mm}

\begin{lemma}
\label{L5.3}
In set ${\cal B},$  $\alpha$ is a continuous function of $(a,\beta,\gamma)$, so there exists a compact ${\cal K}\subset{\cal B}$ not depending on $v>0$ such that $\sup_{{\cal B}}\tilde I^2=\sup_{{\cal K}}\tilde I^2$, and in particular the maximum on ${\cal B}$ is attained.
\end{lemma}
\proof
Let us consider $\sup_{{\cal B}}\tilde I^2$. Notice that $C$ defined by (\ref{defA}) depends on $v$. Since $\tilde J^2(a,\alpha,\beta,m)=-1$ on ${\cal B}$, we have
$$\alpha=\alpha(a,\beta,m):=C(\beta)-(1-\beta)\delta\psi(a)+e^*(\beta)+m-\frac{g}{g'}( m),$$
 therefore we have:

\begin{align}
&\hspace*{-0.5cm} \sup_{(a,\alpha,\beta,m)\in{\cal B}}\tilde I^2(a,\alpha,\beta,m)\nonumber
\\
&=\sup_{\{(a,\beta,m)\in[0,\bar a]\times(-\infty,\bar\beta]\times\mathbb{R}_+;~\alpha(a,\beta,m)\geq 0\}}F(a)-e^{-u\delta\psi(a)} \tilde B(\beta, m) 
 \nonumber\\
&=\sup_{\{(a,\beta,m)\in[0,\bar a]\times(-\infty,\bar\beta]\times\mathbb{R}_+;~\alpha(a,\beta,m)\geq 0\}} F(a)-e^{u(Id-b^m)(m)}e^{-u\delta\psi(a)}\tilde k(\beta) \label{expressionIenbetagamma} 
 \end{align}
with
$$\tilde k(\beta)=e^{u(Id-b^e)\circ e^*(\beta)}    e^{uC(\beta)}\mathbb{E}\left[e^{u\beta\big(\mu~-\varphi\circ e^*(\beta)\big)}\right].$$
In (\ref{expressionIenbetagamma}), we have $F(a)$ minus the product of three positive functions, the first one depending only on $m$, the second one only on $a$ and the third one, $\tilde k$, only on $\beta$. By \eqref{eqm}, 
the first function $e^{u(m-b^m(m))}$ goes to $\infty$ when $m \to\infty$  so it is bounded from below by a positive constant. Since $\delta^+\psi$ is bounded form above, 
the second one is also bounded from below by a positive constant. 
\\
As in  Lemma \ref{L4.3}, we prove that $\tilde k(\beta)\mapsto \infty$ uniformly in $v>0$ when $\beta\mapsto-\infty.$ Since $(0,0,1,0)\in{\cal B}$ for any $v>0$, we deduce the existence of $\bar{m}\in[0,+\infty)$ and $\underline\beta\in (-\infty,\bar\beta]$ not depending on $v>0$ such that $\sup_{\cal B}\tilde I^2=\sup_{{\cal B}\cap\{\beta\in[\underline\beta,\bar\beta],m\in[0,\bar{m}]}\tilde I^2$.To conclude the proof, we remark that
\begin{align*}
   &\sup_{a\in[0,\bar a],\beta\in[\underline \beta,\bar\beta],m\in[0,\bar m]}\tilde I^2(a,\alpha,\beta,m)\leq \sup_{a\in[0,\bar a]}F(a)\\
 &-e^{u\alpha}\inf_{a\in[0,\bar a],\beta\in[\underline \beta,\bar\beta],m\in[0,\bar m]\}}\mathbb{E}\left[e^{-u\big(b^m(m)+b^e(e^*(\beta))-\beta(\mu~-\varphi(e^*(\beta))-\delta \psi(a))-
 \frac{ g}{g'}(m)\big)}\right]
\end{align*}
tends to $-\infty$ as $\alpha\to\infty$.
\ep

\begin{lemma}\label{L5.4}
Assume that 
there exists a Stackelberg equilibrium $(\hat e, \hat m, \hat a, \hat\alpha, \hat\beta,\hat\gamma)$ with $\hat\alpha>0$. Then  necessarily 
\begin{eqnarray}
\label{eqdegamma}
&\hat\gamma=\frac{1}{g'(\hat m)},~\mbox{where}~\hat m=[(b^m)']^{-1}(1)>0,
\\
\label{eqdebeta}
&\hat\beta<\frac{v}{u+v}~\mbox{and}~\hat\beta~\mbox{is a solution of (\ref{eqbetaStackelberg})}.
\end{eqnarray}
\end{lemma}
\proof
If there exists a Stackelberg equilibrium $(\hat e, \hat m, \hat a, \hat\alpha, \hat\beta,\hat\gamma)$, then the supremum of $\tilde I^2$ on ${\cal A}$ is attained at $(\hat a, \hat\alpha, \hat\beta,\hat\gamma)$. Lemma \ref{L5.2} and the hypothesis $\hat \alpha>0$ yield that $\tilde J^2(\hat a, \hat\alpha, \hat\beta,\hat\gamma)=-1$. Since $(b^m)'(0)=+\infty$, $m\mapsto m-b^m(m)$ is decreasing in a neighborhood of $0$, so that from (\ref{expressionIenbetagamma}), $\hat m>0$. The mapping $(\beta,m)\mapsto \tilde I^2(\hat a,\alpha(\hat a,\beta,m),\beta,m)$ where $\alpha(a,\beta,m)$ is the function introduced at the beginning of the proof of Lemma \ref{L5.3} admits a local maximum at $(\hat\beta,\hat m)$. Therefore the first order conditions are satisfied (notice that $\bar{\beta}$ may be increased), i.e. $\frac{\partial}{\partial \beta} \tilde I^2(\hat a,\alpha(\hat a,\hat\beta,\hat m),\hat\beta,\hat m)=\frac{\partial}{\partial m} \tilde I^2(\hat a,\alpha(\hat a,\hat\beta,\hat m),\hat\beta,\hat m)=0$. 
 The partial derivative with respect to $m$ yields
 $$
\hat m=[(b^m)']^{-1}(1)>0.
$$
 On the other hand, using $e^*(\beta))=(\varphi')^{-1}(\frac{1}{(1-\beta)^+}),$ 
 we compute:
\begin{align*}
C'(\beta)=\begin{cases}
-\frac{\E\left[\mu~e^{v(1-\beta)\mu~}\right]}{\E\left[ e^{v(1-\beta)\mu~}\right]}=\varphi(0)-h(v(1-\beta))~\mbox{if}~\beta\geq 1\\
\varphi\circ e^*(\beta)-(e^*)'(\beta)-h(v(1-\beta))~\mbox{if}~\beta<1.
\end{cases}
\end{align*}
Then:
\begin{align*}
\tilde k'(\beta)=\begin{cases}
\E \left[u\Big(\mu+\varphi(0)~-h(v(1-\beta))\Big)e^{u\big(\beta\mu~+C(\beta)\big)}\right]~\mbox{if}~\beta\geq 1\\
u\E\left[\Big(-\frac{\beta}{1-\beta}(e^*)'(\beta)-(b^e)'\circ e^*(\beta)(e^*)'(\beta)+\mu~-h(v(1-\beta))\Big)
\right.\\
\left.
\qquad\qquad e^{u\big((Id-b^e)\circ e^*(\beta)+\beta(\mu~-\varphi\circ e^*(\beta))+C(\beta)\big)}\right]~\mbox{if}~\beta<1.
\end{cases}
\end{align*}
As seen in Lemma \ref{lem:h}, for $\beta\geq 1$,   $h(u\beta)-h(v(1-\beta))>0$, therefore $\tilde k'(\beta)>0$ for $\beta\geq 1$. As a consequence, the equation $\frac{\partial}{\partial\beta} \tilde I^2=\tilde k'(\beta)=0$ implies that $\hat\beta$ solves (\ref{eqbetaStackelberg}). By Proposition \ref{eq5.5}, $\hat\beta<\frac{u}{u+v}$.
\ep


\begin{lemma}
\label{L5.5}
The supremum of $(a,\alpha,\beta,m)\mapsto \tilde I^2(a,\alpha,\beta,m)$ on ${\cal C}$ defined by (\ref{calC}) is attained.
\end{lemma}
\proof
Like in the proof of Lemma \ref{L4.5}, there exists $\bar m\geq m_0$ not depending on $v>0$ such that $\tilde I^2$ is decreasing with respect to $m$ on $[\bar m,+\infty)$. 
  Thus in case $m\geq\bar m,$ the optimum has to  bind the constraint.
 Therefore, writing
$${\cal C}_1:=\{(a,0,\beta,m)\in{\cal C}: m\in[0,\bar m]\},$$
$$ {\cal C}_2:=\{(a,0,\beta,m)\in{\cal C}:m\geq\bar m~\mbox{and}~
 \tilde J^2(a,0,\beta,m)= -1 \},$$
the supremum of $\tilde I^2$ over ${\cal C}$ is attained iff the supremum of $\tilde I^2$ over ${\cal C}_1\cup{\cal C}_2$ is attained.\\

(i)
We prove that ${\cal C}_1$ is compact. The condition $-\tilde J^2(a,0,\beta,m)\leq 1$ is equivalent to
$C(\beta)-(1-\beta)\delta\psi(a)+e^*(\beta)+m-\frac{g}{g'}(m)\leq 0$. 
For $\beta<1$, by \eqref{minoc}, it implies $e^*(\beta)+m-\frac{g}{g'}(m)\leq 0$. Since $\lim_{\beta\to-\infty}e^*(\beta)=+\infty$, while $m\in[0,\bar m]$, this implies the existence of $\underline{\beta}\in (-\infty,\bar\beta]$ such that ${\cal C}_1\subset [0,\bar a]\times\{0\}\times[\underline\beta,\bar\beta]\times[0,\bar m]$ where $\bar m$ and $\underline\beta$ do not depend on $v$. Since 
 ${\cal C}_1$ is closed by continuity of $\tilde J^2$, it is compact and therefore the supremum over ${\cal C}_1$ is attained.

\vspace{5mm}

(ii)
On ${\cal C}_2$, we have $C(\beta)-(1-\beta)\delta\psi(a)+e^*(\beta)=\frac{g(m)}{g'(m)}-m=G(m)$ with $m\geq \bar{m}\geq m_0$ and therefore, by Lemma \ref{Lemmalimginfty}, $m=H(a,\beta):=G^{-1}\big(C(\beta)-(1-\beta)\delta\psi(a)+e^*(\beta)\big)$. We thus have $\sup_{{\cal C}_2}\tilde I^2=\sup_{{\cal C}_2\cap \{a\in [0,\bar{a}],\beta\in(-\infty,\bar\beta]\}}\tilde I^2\big(a,0,\beta,H(a,\beta)\big)$. 
\\
Adapting the proof of Lemma \ref{L4.5}, we check that when $\beta\to-\infty$, $\tilde I^2(a,0,\beta,H(a,\beta))\to-\infty$ uniformly in $(a,v)\in[0,\bar{a}]\times (0,+\infty)$.\\Since $(0,0,1,G^{-1}(0))\in{\cal C}$ for any $v>0$ and $\tilde I^2(0,0,1,G^{-1}(0))$ does not depend on $v>0$, one deduces the existence of $\underline\beta\in (-\infty,\bar\beta]$ not depending on $v>0$ such that if $\sup_{\cal C}\tilde I^2=\sup_{{\cal C}_2}\tilde I^2$ then $\sup_{\cal C}\tilde I^2=\sup_{{\cal C}_2\cap\{\beta\in[\underline\beta,\bar\beta]\}}\tilde I^2$. Now 
 \begin{align*}
   &\sup_{a\in[0,\bar{a}],\beta\in[\underline{\beta},\bar\beta]}\tilde I^2(a,0,\beta,m)\leq\sup_{a\in[0,\bar{a}]}F(a)\\&-e^{u\left(\frac{g}{g'}(m)-b^m(m)\right)}\inf_{a\in[0,\bar{a}],\beta\in[\underline{\beta},\bar\beta]}\mathbb{E}\left[e^{-u\big(b^e(e^*(\beta))-\beta(\mu~-\varphi(e^*(\beta))-\delta\psi(a))\big)}\right]
\end{align*}
where, by \eqref{eqm}, the right-hand-side tends to $-\infty$ as $m\to\infty$. Hence there exists $\bar{m}\in[0,+\infty)$ not depending on $v>0$ such that if $\sup_{\cal C}\tilde I^2=\sup_{{\cal C}_2}\tilde I^2$ then $\sup_{\cal C}\tilde I^2=\sup_{{\cal C}_2\cap\{\beta\in[\underline\beta,\bar\beta],m\in[0,\bar{m}]\}}\tilde I^2$.\ep


\vspace{5mm}

\noindent{\bf End of the proof of Theorem \ref{PropStackelberg}.} 
\\
In conclusion, the maximum over ${\cal A}$ is attained, either on ${\cal B}$ or ${\cal C}$, and a Stackelberg equilibrium exists. Moreover if $\alpha>0$, the characterization given in the statement of  the theorem follows from Lemma \ref{L5.4}. The last assertion follows from the end of the proof of Proposition \ref{PropOptimumFirm1}. \ep

\begin{remark}
(i)
Notice that $\tilde I^2$ is not in general concave w.r.t $\beta$ or $\gamma$. Consider for example $b^m(x)=x^{0.8}/0.8$ and $g(x)=x^{0.5}/0.5$.
\\
(ii)
 Assume that equation (\ref{eqbetaStackelberg}) admits a unique solution, namely $\hat\beta$. Then, from Lemma \ref{L5.4},
any 
Stackelberg equilibrium $(\hat e,\hat m,\hat a,\hat\alpha,\hat\beta,\hat\gamma)$ with $\hat\alpha>0$ satisfies : 
$\hat e=(\varphi')^{-1}\Big(\frac{1}{1-\hat\beta}\Big)$, $\hat m=\big[(b^m)'\big]^{-1}(1)$, $\hat\gamma=1/g'(\hat m)$,
\\
 $\hat a\in\arg\max_{a\geq 0} F(a)-e^{-u\delta\psi(a)}\tilde B(\hat\beta,\hat m)$, $\hat\alpha=C(\hat\beta)-(1-\hat\beta)\delta\psi(\hat a)+\hat e+\hat m-\hat\gamma g(\hat m)$, where $C$ and $\tilde B$ are respectively defined by (\ref{defA}) and (\ref{defB}). 
\\
In particular, if $\delta\geq 0$, an adaptation of the end of the proof of Proposition \ref{PropOptimumFirm1} proves that $\hat a$ is unique. 
 \end{remark}

\section{Proofs of Comparison  results}
\label{sec:proof3}

\subsection{Comparison between both Nash equilibria}
\label{comp}

We summarize Nash equilibria in both situations : in Situation 1, the Nash equilibria is characterized by
$$
\hat\beta=\frac{v}{u+v},~\hat e=(\varphi')^{-1}(\frac{u+v}{u}),~\hat a_1=
(\psi')^{-1}\left(\frac{(u+v)(1+r_J)}{\delta^+ u}\right),~\hat m\geq 0, $$
 $$\hat \gamma=\frac{1}{g'(\hat m)}, ~ \hat\alpha_1+\frac{g(\hat m)}{g'(\hat m)}+\frac{u}{u+v}\delta\psi(\hat a_1)=(1+r_J)(D+\hat a_1)+\hat m+\hat e+C( \hat \beta)$$
  which leads to the optimal value for firm $I$:
 $$\hat I^1(\hat m)=b^a(\hat a_1)-
 {B(\hat m)}e^{u[(1+r_J)(D+\hat a_1)-\delta\psi(\hat a_1)]}
 , $$
  where we recall
$${B(\hat m)}= e^{-u(b^m(\hat m)-\hat m+b^{e}(\hat e)-\hat e)}e^{(u+v) C( \hat \beta)},$$
$$e^{v C( \hat \beta)}=\E\left[ \exp(\frac{uv}{u+v}(\mu~-\varphi(\hat e)))\right].$$
 In Situation 2, the Nash equilibria are characterized by
$$
\hat\beta=\frac{v}{u+v},~
\hat a_2\in\arg\max_{a\geq 0} \left(b^a(a)-f[(1+r_I)(D+a)]-e^{-u\delta\psi(a)}B(\hat m)\right), $$
 $$\hat e=(\varphi')^{-1}(\frac{u+v}{u}),~ \hat\alpha_2+\frac{g(\hat m)}{g'(\hat m)}+\frac{u}{u+v}\delta\psi(\hat a_2)=\hat m+\hat e+ C( \hat \beta),~\hat m\geq 0,~ \hat \gamma=\frac{1}{g'(\hat m)},$$
 which leads to the optimal value for firm $I$:
 $$\hat I^2(\hat m)=b^a(\hat a_2)
-f[(1+r_I)(D+\hat a_2)]-{B(\hat m)}e^{-u\delta\psi(\hat a_2)}.$$

\noindent The following proposition gives the monotonicity of the optimal initial effort   $\hat a_1$ in Situation $1$  (respectively $\hat a_2$ in Situation 2) as a function  of the interest rate $r_J$ (respectively $r_I$).
 \begin{proposition}
 \label{prop:ar}
 The application $r_J\mapsto \hat a_1(r_J)$ is non increasing.
 \\
 At least in case $\delta>0,$ the application $r_I\mapsto \hat a_2(r_I)$ is non increasing.
  \end{proposition}
\proof
 The monotonicity property of the function $r_J\mapsto \hat a_1(r_J)$ is a trivial consequence of the definition of $\hat a_1$,
 since by hypothesis the function $\psi'$ is non increasing.
 In case $\delta<0,$ $\hat a_1(r_J)=0.$\\
 The second assertion is a consequence of the characterization of $\hat a_2(r_I)$ as the unique solution 
 of the equation (see the end of the proof of Proposition \ref{PropOptimumFirm1}) $H(\hat a_2(r_I),r_I)=0$ where
 \begin{equation}
 \label{Defa2}
 H(a,r_I)=(b^a)'(a)+u\delta\psi'(a) B(\hat m) e^{-u\delta\psi(a)}-(1+r_I)f'[(1+r_I)(D+a)].
 \end{equation}
 Thus $\frac{d\hat a_2}{dr_I}(r_I)=-\frac{\partial_{r_I}H}{\partial_a H}(\hat a_2(r_I),r_I).$ 
 Notice that 
 $$\partial_{r_I}H=-f'[(1+r_I)(D+a)]-(1+r_I)(D+a)f''[(1+r_I)(D+a)]<0$$
 since $f$ is convex non decreasing. Therefore the sign of $\frac{d\hat a_2}{dr_I}$
 is the one of $ \partial_a H=$
 \begin{equation}
  (b^a)''(a)+u\delta B(\hat m) e^{-u\delta\psi(a)}[\psi''(a)-u\delta(\psi'(a))^2]-(1+r_I)^2f''[(1+r_I)(D+a)]<0\label{decha}
 \end{equation}
 since $b^a$ and $\psi$ are concave and $f$ is convex.
 \ep
 \\
 
\noindent  We now prove  a sufficient  condition under which the best
  situation is the second one (debt issuance), as stated in Theorem 
 \ref{prop:1stCS}.\\
 
\noindent {\bf  Proof of Theorem \ref{prop:1stCS}}
\\
  The key of the proof is the remark that, since in Situation 2 $a$ is a control of firm $I$, the optimal value $\hat I^2(\hat m)$ obtained by this firm in any Nash equilibrium $(\hat\alpha, \hat\beta, \hat\gamma, \hat e,\hat m, \hat a_2(r_I))$ is larger than $I^2(\hat\alpha, \hat\beta, \hat\gamma, \hat e,\hat m,\hat a_1(r_J))$. Hence
\begin{eqnarray*}
&&  \hat I^1(\hat m)-\hat I^2(\hat m)\leq \hat I^1(\hat m)-I^2(\hat\alpha, \hat\beta, \hat\gamma, \hat e,\hat m,\hat a_1(r_J))
 \\
& &=f((1+r_I)(D+\hat a_1(r_J)))- B(\hat m) e^{-u\delta\psi(\hat a_1(r_J))}(e^{u(1+r_J)(D+\hat a_1(r_J))}-1).
\end{eqnarray*}
Condition (\ref{ccl2}) is equivalent to non-positivity of the right-hand side.
 \ep\\

\noindent In the case
$\delta >0$, we now prove  a sufficient  condition under which the best
  situation is the first one  (outsourcing), as stated in Theorem
 \ref{prop:2dCS}.\\

\noindent {\bf  Proof of Theorem \ref{prop:2dCS}}
\\
Let $\hat m$ be the same parameter in both situations.
By a slight abuse of notations, we introduce the function
$$I^1:~a\mapsto b^a(a)-B(\hat m) e^{u\left((1+r_J)(D+a)-\delta\psi(a)\right)}$$
which is such that $\hat I^1(\hat m)=I^1(\hat a_1(r_J))$.
Let us check that (\ref{CSI^1(a1)(a2)}) implies that
$I^1(\hat a_2(r_J))<I^1(\hat a_1(r_J))$. Since $\delta>0$, the second condition in (\ref{CSI^1(a1)(a2)}) implies that $\hat a_1(r_J)>\hat a_2(r_I)$ by the definition of $\hat a_1(r_J)$ and the monotonicity of $\psi'$. In terms of the function $H$ defined by (\ref{Defa2}), the first condition in (\ref{CSI^1(a1)(a2)}) writes $H(\hat a_1(r_J),r_I)<0.$  Since $H(\hat a_2(r_I),r_I)=0$ and the function $a\mapsto H(a,r_I)$ is decreasing by (\ref{decha}), this also implies that $\hat a_1(r_J)>\hat a_2(r_I)$.
\\
Since $a\mapsto e^{u\left((1+r_J)(D+a)-\delta\psi(a)\right)}$ is convex by
composition of the convex function $a\mapsto u\left((1+r_J)(D+a)-\delta\psi(a)\right)$
 with the increasing and convex exponential function, $I^1$ is concave from 
 the  concavity of $b^a$. Now $(I^1)'(\hat a_1(r_J))=$
$${b_a}'\left((\psi')^{-1}(\frac{(u+v)(1+r_J)}{\delta u})\right)+
 B(\hat m) e^{u\left((1+r_J)(D+a)-\delta\psi(a)\right)}u(1+r_J)>0.$$
From the concavity of $I^1$ and the inequality $\hat a_1(r_J)>\hat a_2(r_I)$, we deduce that $I^1(\hat a_2(r_J))<I^1(\hat a_1(r_J))=\hat I^1(\hat m)$. As a consequence,
\begin{eqnarray*}
 &&  \hat I^1(\hat m)-\hat I^2(\hat m)>I^1(\hat a_2(r_J))-\hat I^2(\hat m)
   \\
   &&=f((1+r_I)(D+\hat a_2(r_I)))- B(\hat m) e^{-u\delta\psi(\hat a_2(r_I))}(e^{u(1+r_J)(D+\hat a_2(r_I))}-1).
\end{eqnarray*}
Condition (\ref{ccl1}) is equivalent to the non-negativity of the right-hand side, which concludes the proof.
 \ep
 \\

\section{Proofs of  incomplete information results}
\label{proofIncInf}
In the incomplete information framework, the firms do not have a perfect knowledge of the preferences of the other firm. More precisely, we still assume that the firms' utility functions are $U(x)=-e^{-ux}$ and $V(x)=-e^{-vx}$ respectively, but firm $I$ perceives $v$ as a $(0,+\infty)$-valued random variable ${\cal V}$ independent of $\mu$ with a known distribution. The penalty that firm $I$ gets if firm $J$ does not accept the contract is denoted by $p\in\mathbb{R}\cup\{+\infty\}$.

\subsection{Stackelberg equilibrium in incomplete information,   firm $I$ is leader}
(cf. (\ref{eqinfoincomplete1}) and  (\ref{eqinfoincomplete2})) 
 \begin{equation}
 \label{**}
u^i_I:=-p\vee\sup_{c^i}\{I^i(c^i)\pr\left({\cal A}^i(c^i)\right)-p(1-\pr\left({\cal A}^i(c^i)\right))\},
 \end{equation}
where $c^i$ is the control of firm I in Situation $i$ :
$c^1= (\alpha,\beta,\gamma)$ and $c^2= (a, \alpha,\beta,\gamma)$ and ${\cal A}^i(c^i)$ is the event ``firm $J$ accepts the contract''. 
Firm $J$ accepts the contract if and only if $- \tilde J^i({\cal V},c^i)\leq 1,$ 
therefore ${\cal A}^i(c^i)=\{\tilde J^i({\cal V},c^i)\leq 1\}$ where $\tilde J^1$ and $\tilde J^2$ are respectively defined in (\ref{deftJv1}) and  (\ref{deftJv2}). The next lemma  aims at expliciting this acceptance set.
 The functions $v\mapsto -\tilde{J}^i(v,c^i)$ are infinitely differentiable according to Hypothesis (\ref{hypmu}) and convex since their second order derivative is non-negative. Since $- \tilde J^i(0,c^i)=1$, one deduces that
\begin{lemma}
\label{lemacc}
For  $c^1=(\alpha,\beta,\gamma)\in \mathbb{R}_+\times\mathbb{R}\times\mathbb{R}_+$ (respectively for $c^2=(a,\alpha,\beta,\gamma)\in\mathbb{R}_+\times\mathbb{R}_+\times\mathbb{R}\times\mathbb{R}_+$), $\bar v(c^i):=\sup\{v\geq 0:- \tilde J^i(v,c^i)\leq 1\}$ belongs to $[0,+\infty]$. If $\bar v(c^i)\in[0,+\infty)$, then $\{v\geq 0:- \tilde J^i(v,c^i)\leq 1\}=[0,\bar v(c^i)]$ and $- \tilde J^i(\bar v(c^i),c^i)=1$. If $\bar v(c^i)=+\infty$, then $\{v\geq 0:- \tilde J^i(v,c^i)\leq 1\}=[0,+\infty)$.
\end{lemma}

\noindent Recall the value function of the problem with complete information that firm J's risk aversion is equal to $v$:
 $$u^i(v):=\sup_{\{c^i:- \tilde J^i(v,c^i)\leq 1\}}\tilde I^i(c^i)$$
   and  $ w^i_I=-p\vee \sup_{v>0}\{u^i(v)\pr({\cal V}\leq v)-p\big(1-\pr({\cal V}\leq v)\big)\}$ as defined in Theorem \ref{Propinfoincomplete}.
We now are able to prove Theorem \ref{Propinfoincomplete}.

\vspace{5mm}

\noindent {\bf {Proof of Theorem \ref{Propinfoincomplete}}}

\noindent 
We first prove that $u^i_I\geq w^i_I$. We only need to do so when $w^i_I>-p$ which implies that 
$$w^i_I=\sup_{v>0:u^i(v)\pr({\cal V}\leq v)-p\pr({\cal V}> v)>-p}\left(u^i(v)\pr({\cal V}\leq v)-p\pr({\cal V}>v)\right).$$ 
Let $v>0$ be such that $u^i(v)\pr({\cal V}\leq v)-p\pr({\cal V}> v)>-p$ and $\hat c^i(v)$ be an optimal control for $u^i(v)$, so that $u^i(v)=\tilde I^i(\hat c^i(v))$. 
Since $-\tilde J(v,\hat c^i(v))\leq 1$, by Lemma \ref{lemacc} one has $\{{\cal V}\leq v\}\subset {\cal A}^i(\hat c^i(v))$ and when $p<+\infty$,
\begin{eqnarray*}
&u^i(v)\pr({\cal V}\leq v)-p\pr({\cal V}>v)=(u^i(v)+p)\pr({\cal V}\leq v)-p
\\&\leq (u^i(v)+p)\pr\left({\cal A}^i(\hat c^i(v))\right)-p
\\
 &=\tilde I^i(\hat c^i(v))\pr({\cal A}^i(\hat c^i(v))-p(1-\pr({\cal A}^i(\hat c^i(v))\leq u^i_I.
\end{eqnarray*}
When $p=+\infty$, the left-most side of the above inequalities is still not greater than the right-most side as $1=\pr({\cal V}\leq v)=\pr({\cal A}^i(\hat c^i(v)))$.
\\
Since $v>0$ such that $u^i(v)\pr({\cal V}\leq v)-p\pr({\cal V}>v)>-p$ is arbitrary, we get $u^i_I\geq w^i_I$.\\
Finally we prove that $w^i_I\geq u^i_I$ in case $u^i_I>-p$. Let $\eps\in(0,\frac{u^i_I+p}{2})$ and $c^i$ be an $\eps$-optimal control for $u^i_I$. Since $$\tilde I^i(c^i)\pr\left({\cal A}^i(c^i)\right)-p(1-\pr\left({\cal A}^i(c^i)\right))\geq u^i_I-\eps>\frac{u^i_I-p}{2}\geq-p,$$ 
one has $\pr\left({\cal A}^i(c^i)\right)>0$. Since ${\cal V}>0$
 and ${\cal A}^i(c^i)=\{{\cal V}\leq \bar v(c^i)\}$ for $\bar v(c^i)=\sup\{v\geq 0:- \tilde J^i(v,c^i)\leq 1\}$, one deduces that $\bar v(c^i)\in(0,\infty)\cup\{\infty\}$.\\
For any $v\in (0,+\infty)$ such that $- \tilde J^i(v,c^i)\leq 1$, one has $u^i(v)\geq \tilde I^i(c^i)$. 
\\
$\bullet$ If $\bar v(c^i)\in (0,+\infty)$ then, by Lemma \ref{lemacc}, $- \tilde J^i(\bar v(c^i),c^i)=1$ so that
\begin{eqnarray*}
  & w^i_I\geq u^i(\bar v(c^i))\pr({\cal V}\leq \bar v(c^i))-p\pr({\cal V}>\bar v(c^i))\geq \tilde I^i(c^i)\pr({\cal V}\leq \bar v(c^i))-p\pr({\cal V}>\bar v(c^i))
   \\
   &=\tilde I^i(c^i)\pr({\cal A}^i(c^i))-p(1-\pr({\cal A}^i(c^i))\geq u^i_I-\eps.
\end{eqnarray*}
$\bullet$
If $\bar v(c^i)=+\infty$, then $\pr\left({\cal A}^i(c^i)\right)=1$ and for all $v>0$, $w^i_I\geq \tilde I^i(c^i)\pr({\cal V}\leq v)-p\pr({\cal V}>v)$ and the same conclusion as before holds by taking the limit $v\to\infty$ in this inequality under the assumption that either $p<+\infty$ or $\pr({\cal V}>v)=0$ for $v$ large enough. Since $\eps>0$ is arbitrarily small, we get $w^i_I\geq u^i_I$, which ends the proof.
 \ep
\\

To prove Theorem \ref{Th:stackincomplet}, we need the following properties of the value functions of the problems with complete information that firm J's risk aversion is $v$.
\begin{lemma}
The function $v\mapsto u^i(v)$ is non-increasing and continuous on $\mathbb{R}_+$.
\end{lemma}
\proof We do the proof for Situation 2. The same holds (with the control $(\alpha,\beta,\gamma)$ instead of $(a,\alpha,\beta,\gamma)$) for Situation 1.
For $v\geq 0$, let ${\cal A}_v=\{(a,\alpha,\beta,\gamma)\in\mathbb{R}_+\times\mathbb{R}_+\times \mathbb{R}\times\mathbb{R}_+:-\tilde J^2(v,a,\alpha,\beta,\gamma)\leq 1\}$. By Lemma \ref{lemacc}, one has ${\cal A}_v\subset{\cal A}_{v'}$ when $v'\leq v$. Therefore $v\mapsto u^2(v)$ is non-increasing.\\
$\bullet$ Let us check the right-continuity of $u^2$ i.e. that $\liminf_{v'\to v^+}u^2(v')\geq u^2(v)$.
According to Proposition \ref{PropStackelberg}, there exists $(a,\alpha,\beta,\gamma)\in{\cal A}_v$ such that $u^2(v)=\tilde I^2(a,\alpha,\beta,\gamma)$. 
\\
Either $-\tilde J^2(v,a,\alpha,\beta,\gamma)<1$ and by continuity of $v'\mapsto -\tilde J^2(v',a,\alpha,\beta,\gamma)$, $(a,\alpha,\beta,\gamma)$
$\in {\cal A}_{v'}$ for $v'$ close enough to $v$ so that the conclusion holds. 
\\
Or $-\tilde J^2(v,a,\alpha,\beta,\gamma)=1$ so that for $v'>v$, $-\tilde J^2(v',a,\alpha,\beta,\gamma)>1,$ and
\\ $\alpha_{v'}=\alpha+\frac{1}{v'}\ln(-\tilde J^2(v',a,\alpha,\beta,\gamma))>0$  is such that $- \tilde J^2(v',a,\alpha_{v'},\beta,\gamma))=1$ and 
$\lim_{v'\to v^+}\tilde I^2(a,\alpha_{v'},\beta,\gamma)=\tilde I^2(a,\alpha,\beta,\gamma)$.\\
$\bullet$ For the left-continuity, we consider a sequence $(v_n)_n$ of positive numbers increasing to a finite limit $v_\infty$.  According to Theorem \ref{PropStackelberg}, there exists $(a_n,\alpha_n,\beta_n,\gamma_n)\in{\cal A}_{v_n}$ such that $u^2(v_n)=\tilde I^2(a_n,\alpha_n,\beta_n,\gamma_n)$.
By Lemma \ref{L5.3} and the proof of Lemma \ref{L5.5}, $(a_n,\alpha_n,\beta_n,\gamma_n)$ stays in a compact subset of $\mathbb{R}_+\times\mathbb{R}_+\times\mathbb{R}\times\mathbb{R}_+$ so one may extract a subsequence that we still index by $n$ for simplicity such that $(a_n,\alpha_n,\beta_n,\gamma_n)$ tends to $(a_\infty,\alpha_\infty,\beta_\infty,\gamma_\infty)$. By continuity of $\tilde I^2$ and $\tilde J^2$, one has $\tilde J^2(v_\infty,a_\infty,\alpha_\infty,\beta_\infty,\gamma_\infty)=\lim_{n\to\infty}\tilde J^2(v_n,a_n,\alpha_n,\beta_n,\gamma_n)$ so that
 $(a_\infty,\alpha_\infty,\beta_\infty,\gamma_\infty)\in{\cal A}_{v_\infty}$ and therefore
$\lim_{n\to\infty}\tilde I^2(a_n,\alpha_n,\beta_n,\gamma_n)=\tilde I^2(a_\infty,\alpha_\infty,\beta_\infty,\gamma_\infty)\leq u^2(v_\infty)$. 
\\
With the monotonicity of $u^2$, we conclude that this function is continuous.
\ep
\\

\noindent We prove the existence of a Stackelberg equilibrium with incomplete information, firm I leader.

\noindent {\bf {Proof of Theorem \ref{Th:stackincomplet}}} 
Let $v_0:=\inf\{v>0:\pr({\cal V}\leq v)>0\}$. 
\begin{description}
   \item[If $\lim_{v\to v_0^+}u^i(v)\leq -p$,], then $v^i_I=w^i_I=-p$.
\item[If not,] $\lim_{v\to v_0^+}u^i(v)>-p$ and we assume that 
$$v_1:=\sup\{v>0:\pr({\cal V}>v)>0\}<+\infty.$$
\begin{description}
   \item[If $p=+\infty$,] then the optimization problem (\ref{Pbequivalent})  clearly admits the solution $v^\star=v_1$.
\item[If $p<+\infty$,] then there exists $v>v_0$ close enough to $v_0$ such that
$$(u^i(v)+p)\pr({\cal V}\leq v)>0.$$Let us deduce existence of a solution to the optimization problem (\ref{Pbequivalent}). Since $u^1$ (resp. $u^2$) is bounded from above by $\sup_{a\in\mathbb{R}_+}b^a(a)<+\infty$ (resp. $\sup_{a\in\mathbb{R}_+}F(a)<+\infty$) and ${\cal V}$ takes its values in $(0,+\infty)$, one has $$\lim_{v\to 0}(u^i(v)+p)\pr({\cal V}\leq v)=0.$$ The function $v\mapsto (u^i(v)+p)\pr({\cal V}\leq v)$ being upper-semicontinuous on the closed set $\{v\in[\varepsilon,v_1]:u_i(v)+p\geq 0\}$ for each $\varepsilon>0$, we conclude that the optimization problem (\ref{Pbequivalent}) has a solution $v^\star\in(0,v_1]$ if $v_0=0$ and in $[v_0,v_1]$ otherwise. 
 \\
 Moreover $u^i(v^\star)>-p$.
Let $\hat c^i(v^\star)$ be an optimal control for $u^i(v^\star)$ such that $u^i(v^\star)=\tilde I^i(\hat c^i(v^\star))$. Since $-\tilde J^i(v^\star,\hat c^i(v^\star))\leq 1$, by Lemma \ref{lemacc} one has 
\\
$\{{\cal V}\leq v^\star\}\subset {\cal A}^i(\hat c^i(v^\star))$ and 
\begin{eqnarray*}
w^i_I&=&u^i(v^\star)\pr({\cal V}\leq v^\star)-p\pr({\cal V}>v^\star)=(u^i(v^\star)+p)\pr({\cal V}\leq v^\star)-p
\\
&\leq &(u^i(v^\star)+p)\pr\left({\cal A}^i(\hat c^i(v^\star))\right)-p
\\
 &=&\tilde I^i(\hat c^i(v^\star))\pr({\cal A}^i(\hat c^i(v^\star)))-p(1-\pr({\cal A}^i(\hat c^i(v^\star)))
 \\
&\leq &u^i_I.
\end{eqnarray*}
With Theorem \ref{Propinfoincomplete}, we conclude that $(\hat c^i(v^\star))$ solves problem (\ref{eqinfoincomplete1})-(\ref{eqinfoincomplete2}).
\end{description}
\end{description}
\ep

\subsection{Nash equilibrium in incomplete information}
We consider both situations, the proofs are quite similar. Let 
$I(a,\alpha,\beta,\gamma,e,m):=$
$$\E \left[ b^a(a) - e^{-u\big(b^m(m)+b^e(e)-\alpha-\beta(\mu~-\varphi(e)-\delta \psi(a))-\gamma g(m)\big)} \right]$$
respectively $I(a,\alpha,\beta,\gamma,e,m):=$
$$\E \left[ F(a) - e^{-u\big(b^m(m)+b^e(e)-\alpha-\beta(\mu~-\varphi(e)-\delta \psi(a))-\gamma g(m)\big)}\right] ,$$
and $J(v,a,\alpha,\beta,\gamma,e,m):=$
$$-\E\left[ e^{-v\big(\alpha+(\beta-1)(\mu~-\varphi(e)-\delta\psi(a))-e+\gamma{g}(m)-m -(1+r_J)(D+ a)\big)}\right]$$
respectively
$$J(v,a,\alpha,\beta,\gamma,e,m):=-\E\left[ e^{-v\big(\alpha+(\beta-1)(\mu~-\varphi(e)-\delta\psi(a))-e+\gamma{g}(m)-m \big)}\right].$$
For firm $I$, given the controls $(a,e,m)$ of firm $J$
(respectively $(e,m)$), the problem is to find $(\alpha,\beta,\gamma)$ 
(respectively $(a,\alpha,\beta,\gamma)$ ) maximizing
$$I(a,\alpha,\beta,\gamma,e,m)\pr(J({\cal V},a,\alpha,\beta,\gamma,e,m)\geq - 1)-p\pr(J({\cal V},a,\alpha,\beta,\gamma,e,m)<-1).$$
As in Lemma \ref{lemacc}, we have
\begin{lemma}
\label{lemacc'}
For $c=(a,\alpha,\beta,\gamma,e,m)\in\mathbb{R}_+\times\mathbb{R}_+\times\mathbb{R}\times\mathbb{R}_+\times\mathbb{R}_+\times\mathbb{R}_+$, $\bar v(c):=\sup\{v\geq 0:J(v,c)\geq -1\}$ belongs to $[0,+\infty]$. If $\bar v(c)\in[0,+\infty)$, then $\{v\geq 0:J(v,c)\geq -1\}=[0,\bar v(c)]$ and $J(\bar v(c),c)=-1$. If $\bar v(c)=+\infty$, then $\{v\geq 0: J(v,c)\geq -1\}=[0,+\infty)$.
\end{lemma}

\noindent {\bf{Proof of Theorem \ref{prop:Nashincomplet}}} 
Assume the  existence of a Nash equilibrium 
$\hat{c}=(\hat{a},\hat{\alpha},\hat{\beta},\hat{\gamma},\hat{e},\hat{m})$ such that the value for firm $I$ is greater than $-p$. This implies that $I(\hat{c})>-p$ and $0<\pr(-J({\cal V},\hat{c})\leq 1)$.
\\
Since by Lemma \ref{lemacc'}, $\pr(-J({\cal V},\hat{c})\leq 1)=\pr({\cal V}\leq \bar v(\hat c))$, one has $\bar v(\hat{c})>0$.\\
We detail below the proof in Situation 1. The one in Situation 2 follows the same scheme, replacing the sets of control parameters $(\alpha,\beta,\gamma)$ and $(a,e,m)$ respectively by $(a,\alpha,\beta,\gamma)$ and $(e,m)$.\\
$\bullet$
Assume that $\bar v(\hat c)<+\infty$ and let $v\in (0,\bar v(\hat c)]$ be such that $\pr({\cal V}\in (v,\bar v(\hat c)])=0$. 
\\
For $(\alpha,\beta,\gamma)$  such that $-J(v,\hat a,\alpha,\beta,\gamma,\hat e,\hat m)\leq 1$:

either
$I(\hat a,\alpha,\beta,\gamma,\hat e,\hat m)\leq -p$ and therefore $I(\hat a,\alpha,\beta,\gamma,\hat e,\hat m)<I(\hat{c}),$ 

or $I(\hat a,\alpha,\beta,\gamma,\hat e,\hat m)>-p$ and  since 
\\
$\pr(-J({\cal V},\hat a,\alpha,\beta,\gamma,\hat e,\hat m)\leq 1)\geq \pr({\cal V}\leq v)=\pr({\cal V}\leq \bar v(\hat c))=\pr(-J({\cal V},\hat{c})\leq 1)$, then:
\begin{eqnarray*}
&&I(\hat a,\alpha,\beta,\gamma,\hat e,\hat m)\pr({\cal V}\leq \bar v(\hat c))-p\pr({\cal V}>\bar v(\hat c))
\\
&\leq& 
   I(\hat a,\alpha,\beta,\gamma,\hat e,\hat m)\pr(-J({\cal V},\hat a,\alpha,\beta,\gamma,\hat e,\hat m)\leq 1)-p\pr(-J({\cal V},\hat a,\alpha,\beta,\gamma,\hat e,\hat m)>1)
   \\
&\leq& I(\hat{c})\pr(-J({\cal V},\hat{c})\leq 1)-p\pr(-J({\cal V},\hat{c})>1)=I(\hat{c})\pr({\cal V}\leq \bar v(\hat c))-p\pr({\cal V}>\bar v(\hat c)),
\end{eqnarray*}
where the last inequality follows from the fact that $\hat{c}$ is a Nash equilibrium for the problem with incomplete information (cf. (\ref{**})).
\\
This implies that $I(\hat a,\alpha,\beta,\gamma,\hat e,\hat m)\leq I(\hat{c})$. Since by Lemma \ref{lemacc'}, $J(v,\hat{c})\geq -1$, we deduce that $\hat{c}$ is a Nash equilibrium for the problem with complete information and risk aversion $v$ for firm $J$. By Theorem 
\ref{propnash2}, we deduce that $\hat{\beta}=\frac{v}{u+v}$ so that the only $v\in (0,\bar v(\hat c)]$ such that $\pr({\cal V}\in (v,\bar v(\hat c)])=0$ is $\bar v(\hat c)=\hat{v}$.\\
$\bullet$
The same line of reasoning permits to conclude that in case $\bar v(\hat c)=+\infty$, $\forall v\in (0,+\infty)$, $\pr({\cal V}\in (v,+\infty))>0$.
\ep


\end{document}